\documentclass[12pt]{article}
\usepackage{graphicx}
\usepackage{booktabs}
\usepackage[a4paper, total={5.5in, 9in}]{geometry}
\usepackage{array}
\usepackage{caption}
\usepackage{float}
\usepackage{multirow}
\usepackage{amsmath}
\usepackage{amssymb}
\usepackage{setspace}
\usepackage{pifont}
\usepackage{breqn}
\usepackage{multicol}
\usepackage{color}
\usepackage{float}
\usepackage[round]{natbib}
\usepackage{hyperref}
\usepackage{geometry} 
\usepackage{algorithm}
\usepackage{algorithmic}
\usepackage{tikz}
\usepackage{pgfgantt}
\usepackage{wrapfig}
\usepackage{adjustbox}
\usepackage{amsmath}
\usepackage{pst-node}
\usepackage{pstricks}
\usepackage{subcaption}
\usepackage{amsmath} 
\usepackage{tikz-qtree}
\usepackage{smartdiagram}
\usepackage{tikz}
\usetikzlibrary{shapes,arrows,calc,arrows.meta}
\usepackage{amsmath,bm,times}
\usepackage{times}
\usepackage{fancyhdr,graphicx,amsmath,amssymb}
\usepackage[title]{appendix}
\usepackage{caption}
\usepackage[font=small,labelfont=bf]{caption}
\usepackage{pdfpages}
\usepackage{booktabs}
\usepackage{pdflscape}
\usepackage{float}
\usepackage[english]{babel}
\usepackage{url}
\usepackage{amsmath}
\usepackage{amssymb}
\usepackage{xcolor}
\usepackage{colortbl}
\usepackage{tabulary}
\usepackage{threeparttable}
\usepackage[utf8]{inputenc}
\usepackage{fancyhdr}
\usepackage[title]{appendix} 
\title{{\Huge  A lexicographic maximin approach to the selective assessment routing problem}}
\author{\large Mohammadmehdi Hakimifar, Vera Hemmelmayr, Fabien Tricoire \vspace{3ex} \\  Institute for Transport and Logistics Management,\\ WU (Vienna University of Economics and
Business), \\Welthandelsplatz 1, 1020 Vienna, Austria\\\href{mhakimif@wu.ac.at}{\{mhakimif, }\href{vera.hemmelmayr@wu.ac.at}vera.hemmelmayr, \href{fabien.tricoire@wu.ac.at}{fabien.tricoire\}@wu.ac.at}}
\date{\small{Keywords: Equity in Humanitarian logistics, Humanitarian needs assessment,\\  Multi-directional local search,  Team orienteering}}

\begin{document}

\vspace{\fill}
\maketitle

\begin{abstract}
	
	 Max-min approaches have been widely applied to address equity as an essential consideration in humanitarian operations.   These approaches, however, have a significant drawback of being neutral when it comes to solutions with the same minimum values.  These equivalent solutions, from a max-min point of view,  might be significantly different. We address this problem using the lexicographic maximin approach, a refinement of the classic max-min approach.  We apply this approach in the rapid needs assessment process, which is carried out immediately after the onset of a disaster, to investigate the disaster's impact on the affected community groups through field visits. We construct routes for an assessment plan to cover community groups, each carrying a distinct characteristic, such that  the vector of coverage ratios are maximized. We define the leximin selective assessment problem, which considers the bi-objective optimization of total assessment time and coverage ratio vector maximization. We solve the bi-objective problem by a heuristic approach based on the multi-directional local search framework. \\
	 
\end{abstract}

\pagenumbering{gobble}
\newpage
\pagenumbering{arabic}

\section{Introduction}
\label{intro}

The max-min (or min-max) approach is one of the widely applied measures to address equity in humanitarian relief operations, such as maximizing the minimum satisfied demand in relief distribution \citep[e.g.,][]{tzeng2007multi,vitoriano2011multi,ransikarbum2016goal,cao2016challenges} or minimizing the maximum arrival time of critical supplies \citep[e.g.,][]{campbell2008routing,nolz2010bi}. A big flaw of this approach is that all solutions with the same max-min value, which might lead to significantly different solutions,  are considered equivalent. This paper aims to refine this approach by applying the lexicographic max-min approach, which has been applied in various supply chain problems such as facility location and network optimization  \citep{ogryczak2014fair}. We focus on rapid needs assessment operations, which in contrast to relief distribution problems, have received little attention in humanitarian logistics literature \citep{pamukcu2020multi}.  According to \citet{luis2012disaster}, majority of studies addressing relief distribution problem assume that the needs of impacted populations are already known or can be estimated, and needs assessment operations are not specifically addressed. 

Rapid needs assessment starts immediately after a disaster strikes to quickly evaluate the disaster's impact by having direct observations of affected sites and conducting interviews with impacted community groups \citep{IFRC2008,ACAPS2011,arii2013rapid}.  The assessment teams have experts from various areas ranging from those who are familiar with the local area to specialties in public health, epidemiology, nutrition, logistics, and shelter \citep{ACAPS2011,arii2013rapid}. A successful rapid needs assessment assists humanitarian agencies in effectively satisfying the needs of affected people at the time of great need \citep{ arii2013rapid}. Planning the field visits and accordingly deciding which sites to visit is an  important decision to reach a  successful assessment.   Since time and resources are limited during the rapid needs assessment stage, it is usually not possible to evaluate the entire affected region; therefore, assessment teams take a sample of affected sites to visit using sampling methods \citep{ACAPS2011}. These methods, in general, aim to select a limited number of sites to visit, which will provide the opportunity for assessment teams to observe and compare the conditions of various affected community groups such as; displaced persons, host communities, and returnees \citep{IFRC2008, ACAPS2011}. The goal here is to achieve acceptable coverage of various community groups.
  
\citet{balcik2017site} propose “Selective Assessment Routing Problem” (SARP), a mixed-integer model that simultaneously addresses site selection and routing decisions of the rapid needs assessment teams. \citet{balcik2017site} aims to evaluate the post-disaster conditions of different community groups, each carrying a distinct characteristic. The objective function in the SARP is maximizing the minimum coverage ratio achieved across the community groups (max-min approach).  The coverage ratio of a characteristic is calculated by  the number of times that the characteristic is covered in the assessment plan divided by the total number of sites in the network carrying that characteristic \citep{balcik2017site}. The SARP has route duration constraints, which define the assessment deadline. 

\begin{figure}[H]
\begin{subfigure}{.5\textwidth}
  \centering
  \includegraphics[width=1\linewidth]{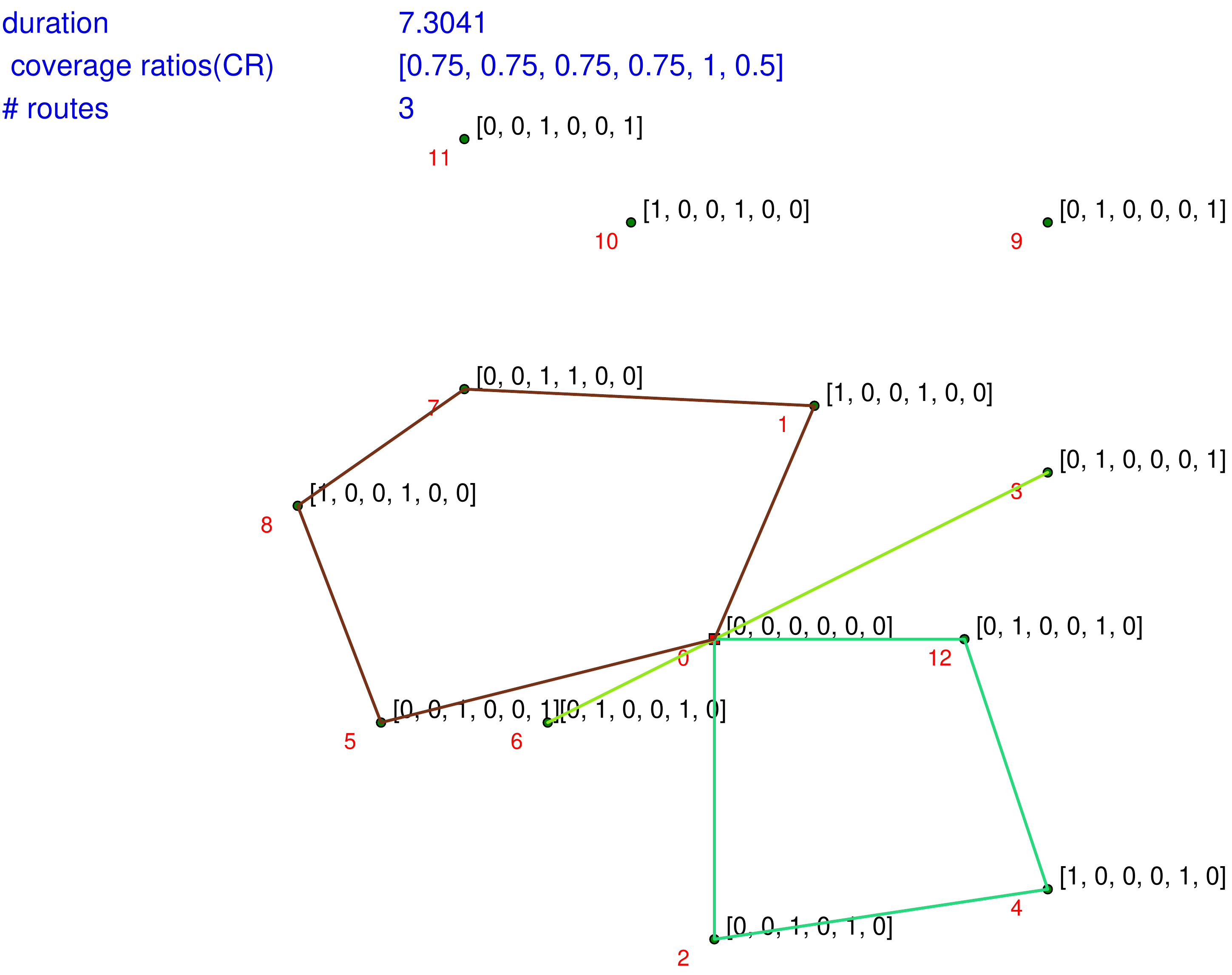}  
  \caption{solution s1}
  \label{fig:sub-first}
\end{subfigure}
\begin{subfigure}{.5\textwidth}
  \centering
  \includegraphics[width=1\linewidth]{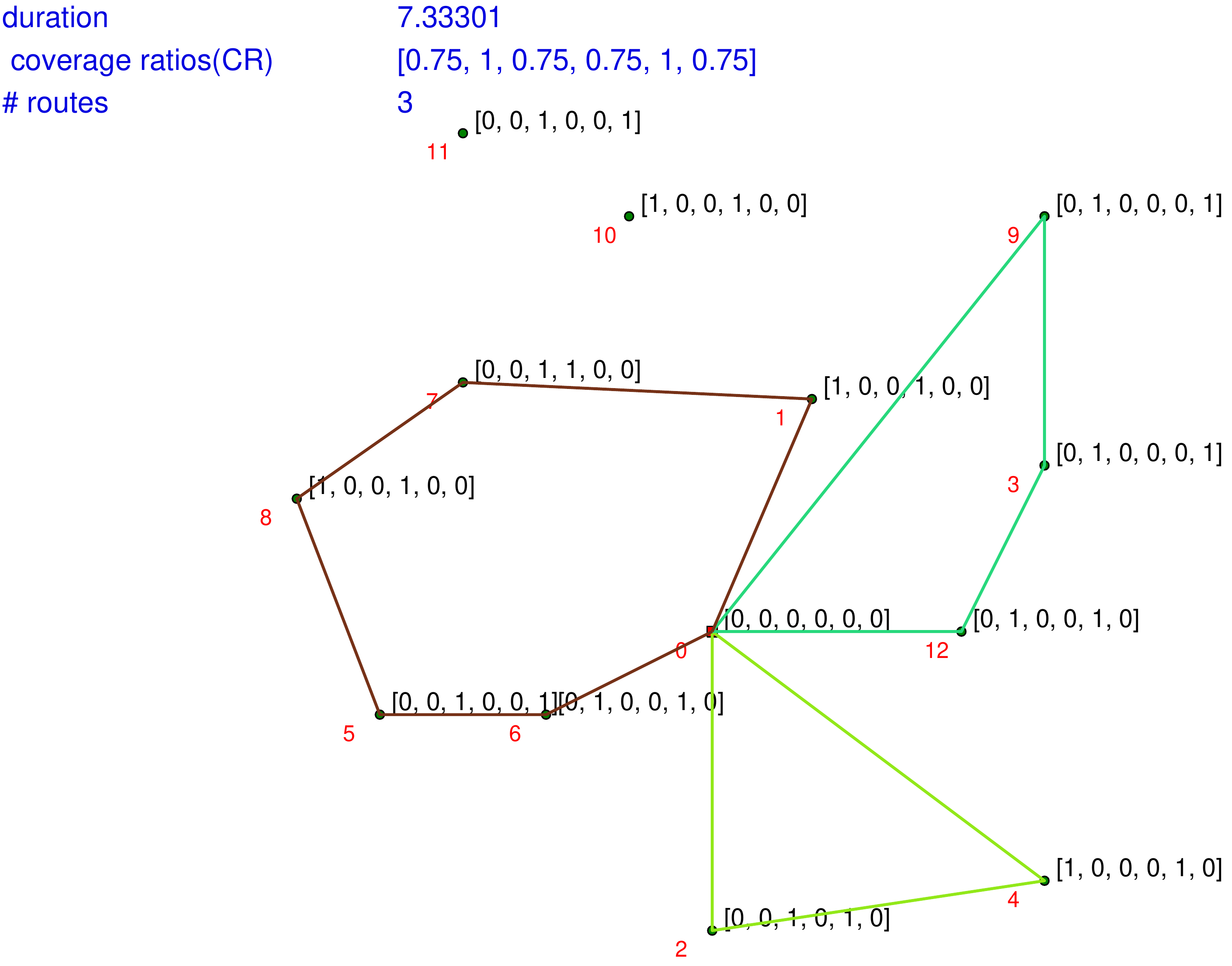}  
\caption{solution s2}
  \label{fig:sub-second}
\end{subfigure}
\begin{center}
\begin{subfigure}{.5\textwidth}
  \centering
  \includegraphics[width=1\linewidth]{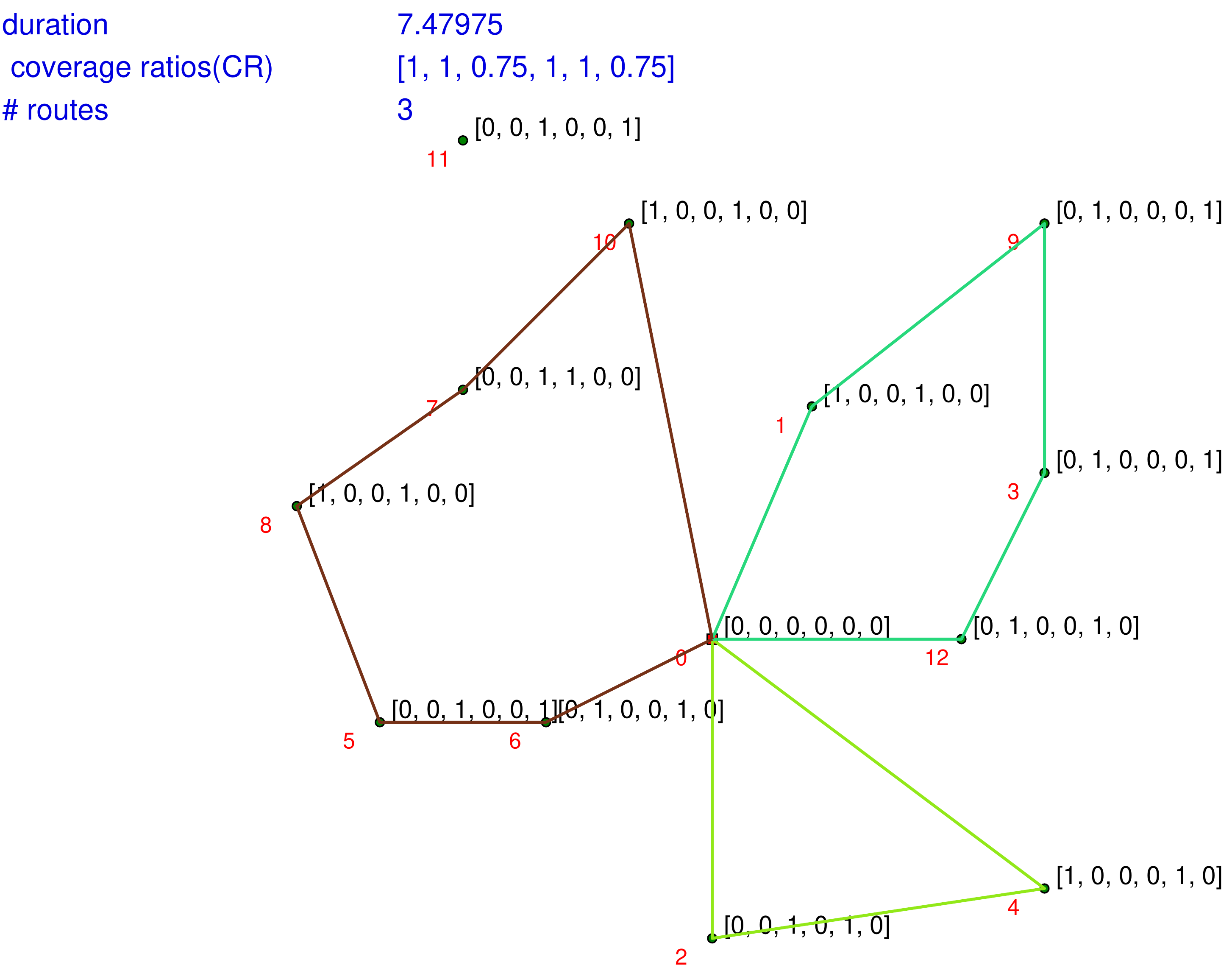}  
\caption{solution s3}
  \label{fig:sub-third}
\end{subfigure}
\end{center}
\caption{Illustrative examples for max-min approach}
\label{maximin}
\end{figure}

 To illustrate the max-min approach and illustrate problems that may arise when not discriminating between solutions with the same minimum coverage ratio in the SARP, we provide an illustrative example in Fig~\ref{maximin}. In this example, there are three solutions ($s_1$, $s_2$, and $s_3$) to a SARP instance with three assessment teams and six characteristics. The array of 0s and 1s at each site indicates which characteristics are present at that site. For example, site number 10 carries characteristics 1 and 4. For solution $s_1$(see Fig~\ref{maximin}-(a)), the minimum coverage ratio is 0.5 whereas in solution $s_2$ (see Fig~ \ref{maximin}-(b)), the minimum coverage ratio is 0.75. Therefore, from a max-min point of view, $s_2$ outperforms $s_1$ since it has a higher minimum coverage ratio. On the other hand, both solutions $s_2$ and $s_3$ (see Fig~ \ref{maximin}-(b) and Fig~ \ref{maximin}-(c)) provide the same minimum coverage value of 0.75, therefore are considered equal from the max-min point of view. However, these two solutions may not be the same from a decision maker's point of view. 
 
In practice, solutions with the same max-min value often have different desirability due to other criteria.  In the case of the SARP, solutions with the same max-min value can provide different coverage for characteristics that are not the least covered characteristic, which is also important: in rapid needs assessment, visiting such extra sites is desirable since each additional visit improves the coverage ratio associated with at least one characteristic \citep{balcik2017site}. Therefore, we believe that introducing the lexicographic maximin approach can be especially relevant in this context. This approach has already been integrated with a few  VRP algorithms \citep{saliba2006heuristics, lehuede2020lexicographic}. 

 During the rapid needs assessment operations, on the one hand, minimizing the total travel time (or having a “faster” assessment) and on the other hand maximizing the reliability level (or having a “better” assessment) are two conflicting but at the same time important goals of decision makers.  As \citet{zissman2014development} indicate, in planning a humanitarian needs assessment for Haiti, the decision makers faced a number of design trade-offs including having a faster or better assessment.“The key challenge was that the data were required urgently: a better assessment could be produced by spending more time designing its elements and verifying its results, but the data were needed faster to inform urgent actions to provide relief to those in need“  \citep[p.13]{zissman2014development}.
In order to capture these two conflicting objectives, we propose a bi-objective optimization problem where we find the set of solutions to the SARP that offer the best trade-off between total route duration (assessment time) and community groups coverage ratio maximization using lexicographic ordering. 

To solve the bi-objective problem, we apply a multi-directional local search method which is a multi-objective optimization framework that generalizes the concept of local search to multiple objectives \citep{tricoire2012multi}. Within this method, we design Adaptive Large Neighborhood Search (ALNS) operators to guide the search towards solutions with higher coverage ratios . We evaluate and analyze the proposed approaches via several computational experiments.

The remainder of the paper is organized as follows. In Section~\ref{sec:lit}, we review the relevant literature. In Section~\ref{lexicoApproach}, we introduce the  leximin approach to the SARP and explain how can we represent the problem in a bi-objective fashion. We present a solution method and computational results in Section~\ref{SolutionMethod} and \ref{Computational}. Finally, we provide  conclusions in Section~\ref{conclusion}.

\section{Related literature}
\label{sec:lit}

Equity has been considered in numerous application areas of
humanitarian supply chains. According to \citet{young1995equity} as cited in \citet{cao2016challenges}, the equity measures based on a different theory of justice can be categorized into three classes: ``(i) Aristotle's proportional equity principle, where the resources are allocated in proportion to the beneficiaries' needs, (ii) classical utilitarianism, where the main aim is to maximize the total welfare of the whole system, and (iii) the Rawlsian maximin or difference
principle, which refers to maximizing the welfare of the worst-off beneficiaries in the system" \citep[34]{cao2016challenges}. In this paper, we focus on the last class of measures within the context of humanitarian operations. In fact, the focus of this paper is not to compare different equity measures as the definition of an appropriate measure depends on the context. 
 
The min-max approach has been applied in numerous relief distribution problems. \citet{tzeng2007multi} construct a relief-distribution model using a multi-objective model for a joint facility location and commodity transportation problem, where the objectives are minimizing
 total costs, minimizing  total travel time, and maximizing the minimal satisfaction during the planning period. \citet{campbell2008routing} focus on two objectives of minimizing the maximum
arrival time of supplies and the average arrival times of supplies. The authors demonstrate the potential impact of using these objective functions on the trade-off between efficiency and equity within the context of the traveling salesperson and vehicle routing
problems. \citet{balcik2008last} address equity by minimizing the maximum percentage of unsatisfied demand over all locations, whereas \citet{ransikarbum2016goal} maximize the minimum percentage of satisfied demand in relief distribution.

Apart from studies that focus on the distribution of relief items,  field visit planning of needs assessment processes has been studied in a few papers in the field of optimization. \citet{huang2013continuous} consider routing of post-disaster assessment teams. In this study, the assessment teams visit all communities in the affected regions, which makes this model appropriate for the detailed assessment stage where  time allows visiting all the impacted sites. This, however, is not the case in the rapid needs assessment stage, where it is usually not possible to visit all the sites due to a shortage of time.  Constructing routes for rapid needs assessment is studied by \citet{balcik2017site}, in which assessment teams need to select a subset of sites to visit and consider the equity among various community groups. \citet{balcik2017site}  presents a  mixed-integer model for the proposed SARP, which simultaneously addresses site selection and routing decisions and supports the rapid needs assessment process that involves the purposive sampling method, a method that only selects those sites that carry certain characteristics.  The objective function is the maximization of the minimum coverage ratio achieved across all community groups.  Minimizing total route duration is used as a secondary objective in a lexicographic fashion.

The study of \citet{balcik2017site} has been extended in various directions. \citet{balcik2020robust} consider the travel time as an uncertain parameter in post-disaster networks and propose a robust optimization approach to address the uncertainty.  \citet{pamukcu2020multi} assume assessment coordinators and specify a target number of observations for each community group in advance, and the objective is to minimize the total duration required to cover all community groups.   \citet{li2020integrated}  present a bi-objective model integrating both the rapid needs assessment stage and the detailed needs assessment stage. \citet{li2020integrated}, similar to \citet{balcik2017site},  maximize the minimum coverage ratio achieved among community groups as the objective of the rapid needs assessment stage.  In the detailed need assessment stage, where the assessment teams need to visit all of the
affected sites, \citet{li2020integrated} minimize the maximum assessment time of all assessment teams. \citet{bruni2020selective} address the post-disaster assessment operations from a customer-centric point of view by considering a service level constraint that guarantees a specified coverage ratio with the objective of minimizing the total latency. \citet{hakimifar2021evaluation} present easy-to-implement heuristic algorithms for site selection and routing decisions of the assessment teams while planning and executing the field visits. They test the performance of proposed heuristic algorithms within a simulation environment and incorporate various uncertain aspects of the post-disaster environment in the field. The authors compare the deterministic setting of the proposed simple heuristics with that of \citet{balcik2017site}.

There are also relevant studies focusing on damage assessment using unmanned aerial vehicles (UAVs) \citep[e.g.,][]{oruc2018post, zhu2019multi,zhu2020optimization,glock2020mission}.  Damage assessment using UAVs studies focus mainly on settings where UAVs' high-quality pictures can meet the assessment purposes, and there is no possibility or necessity to conduct interviews with the community groups, which represent an essential part of needs assessment \citep{hakimifar2021evaluation}. Both SARP and damage assessment with UAVs studies belong to the family of Team Orienteering Problems  \citep{, butt1994heuristic,chao1996team} (TOP).  TOPs address site selection and vehicle routing decisions to maximize the benefits collected from the visited nodes and construct efficient routes. The difference between the TOP and the SARP  lies in the concept of benefits, i.e., the TOP maximizes total collected profits, which are easy to quantify and unique per location, e.g., revenue, whereas the SARP maximizes the minimum coverage ratio of each characteristic over all visited sites. It is not straightforward to quantify the benefits of visiting a particular site without considering the characteristics of the other sites that might be included in the assessment plan  \citep{balcik2017site}. For more information regarding modeling and solution approaches for the TOP, we refer to the surveys conducted by \citet{vansteenwegen2011orienteering}, \citet{archetti2014chapter}, and \citet{gunawan2016orienteering} .

 As an alternative to the classic max-min approach, the lexicographic minimax approach has been applied in various areas of operations research. According to the survey conducted by \citet{ogryczak2014fair}, the lexicographic minimax approach has been applied  in areas such as  facility location \citep{ogryczak1997lexicographic}, supply chain optimization \citep{liu2013multiobjective}, air traffic flow management \citep{bertsimas2008air} and bandwidth allocation and network optimization \citep[e.g.,][]{nace2007lexicographically, ogryczak2005telecommunications}, in order to provide equitable or  fair  solutions.  The lexicographic minimax approach has also been integrated in vehicle routing models. \citet{saliba2006heuristics}  presents heuristics for the VRP with lexicographic max-order objective in a  single-objective fashion.  The lexicographic max-order objective is used  to minimize lexicographically the sorted route lengths. \citet{lehuede2020lexicographic} investigate a lexicographic minimax approach to solve a bi-objective vehicle routing problem with route balancing.  To the best of our knowledge, this study is the first one to introduce the lexicographic maximin approach to the humanitarian operation context. Furthermore,  since the SARP belongs to the family of TOPs, we are also extending the literature on TOPs by introducing the lexicographic maximin approach.

Our contribution is twofold: i) we consider the bi-objective optimization of  finding the set of solutions to the SARP that offer the best trade-off between total route duration and coverage ratio vector maximization. The proposed bi-objective problem addresses the two conflicting objectives of faster versus better assessment that decision makers face during the rapid needs assessment processes. ii) we consider the lexicographic maximin approach -- a refinement of the classic max-min approach.  The lexicographic maximin approach differentiates between solutions with the same max-min value. This differentiation provides a range of solutions for decision-makers,  in which more sites are included in the assessment plan.  Visiting such extra sites is beneficial since it increases the coverage ratio of at least one characteristic.

\section{A lexicographic maximin approach to the selective assessment routing problem}
\label{lexicoApproach}

The lexicographic maximin is a refinement of the max-min approach \citep{dubois1996refinements}. 
In the lexicographic maximin we compare two vectors, instead of simply comparing two numbers (max-min value).
For that purpose, we sort the coverage ratios of a solution from lowest to highest and compare it with the sorted values of the other solution. To illustrate this, let us take solutions $s_2$ and $s_3$ from the illustrative example provided in Fig~\ref{maximin}. The vectors of coverage ratios  from $s_2$ and $s_3$  are shown as bar charts in Figure~\ref{lexicoA} and Figure~\ref{lexicoB} respectively. The sorted forms of these vectors are represented in Figure~\ref{lexicoC} and Figure~\ref{lexicoD}. According to the lexicographic maximin approach, solution  $s_3$ dominates $s_2$ as it has a higher coverage ratio in the third position of sorted vectors. 

We use the following definition of the leximin ordering, adapted from the leximin definition of \citet{bouveret2009computing}:
 
 \begin{figure}[]
\begin{subfigure}{.5\textwidth}
  \centering
  \includegraphics[width=.9\linewidth]{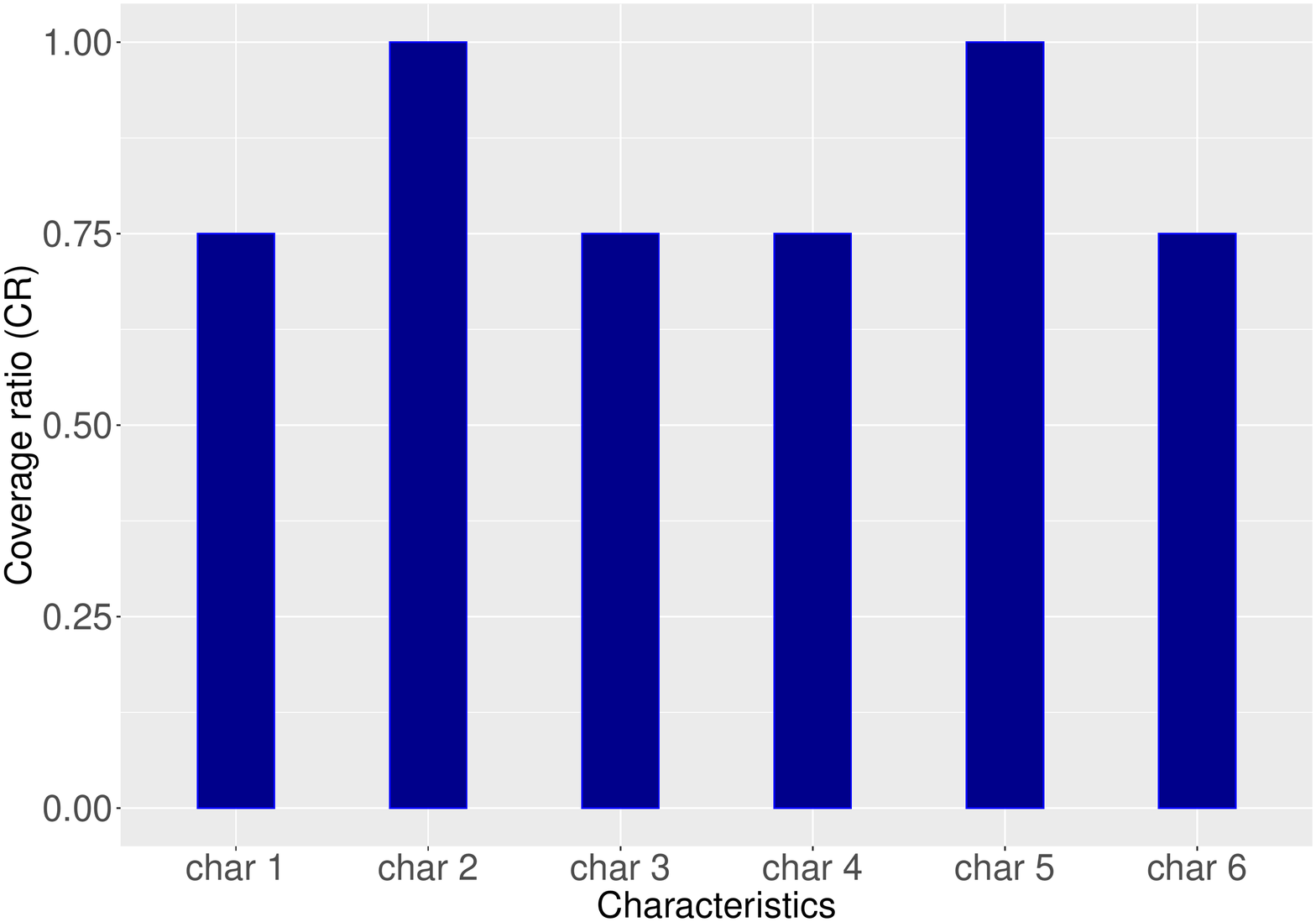}  
  \caption{\scriptsize{s2 coverage ratios: [0.75, 1, 0.75, 0.75, 1, 0.75]}}
  \label{lexicoA}
\end{subfigure}
\begin{subfigure}{.5\textwidth}
  \centering
  \includegraphics[width=.9\linewidth]{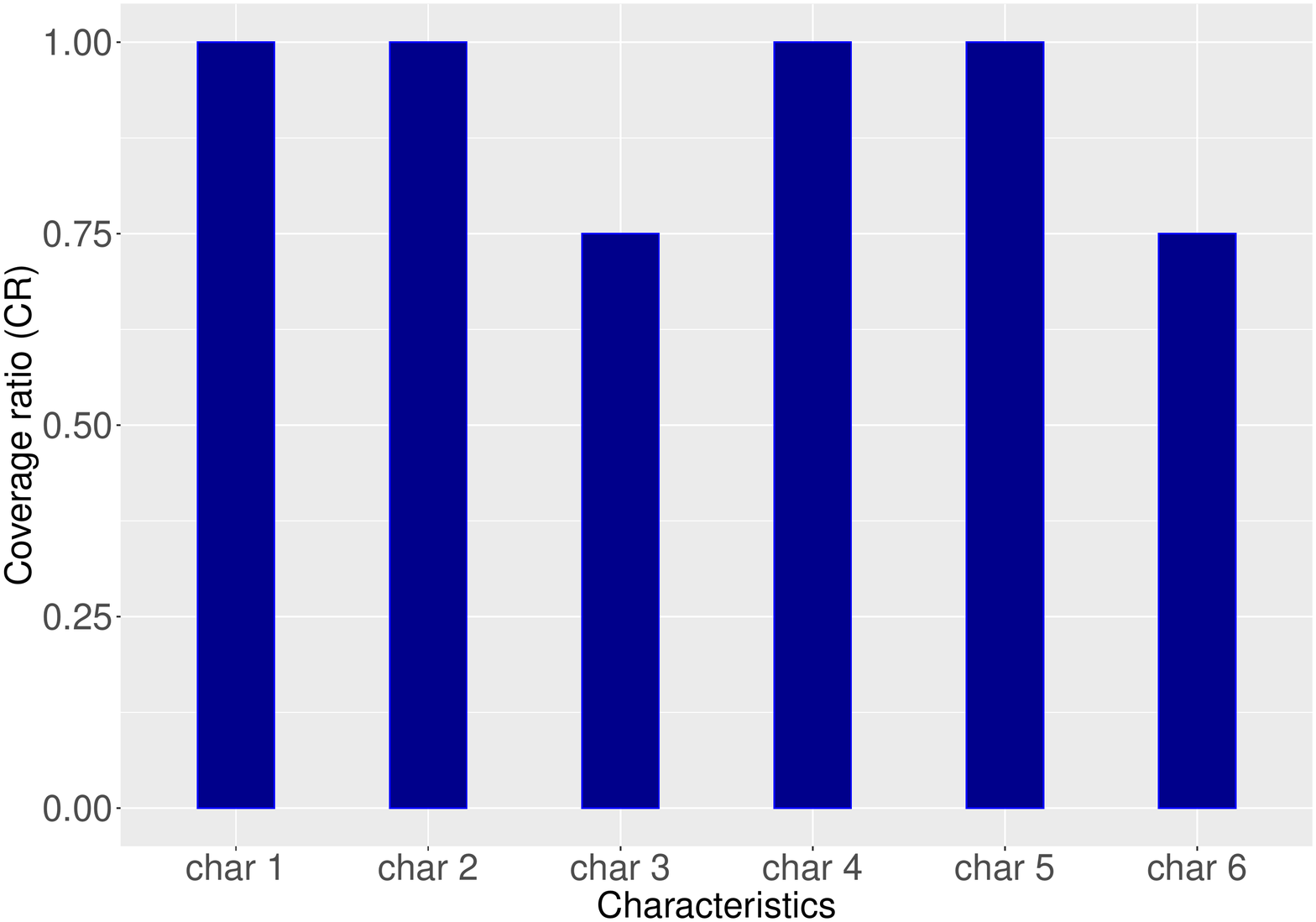}  
\caption{\scriptsize{s3 coverage ratios: [1, 1, 0.75, 1, 1, 0.75]}}
  \label{lexicoB}
\end{subfigure}
\begin{subfigure}{.5\textwidth}
  \centering
  \includegraphics[width=.9\linewidth]{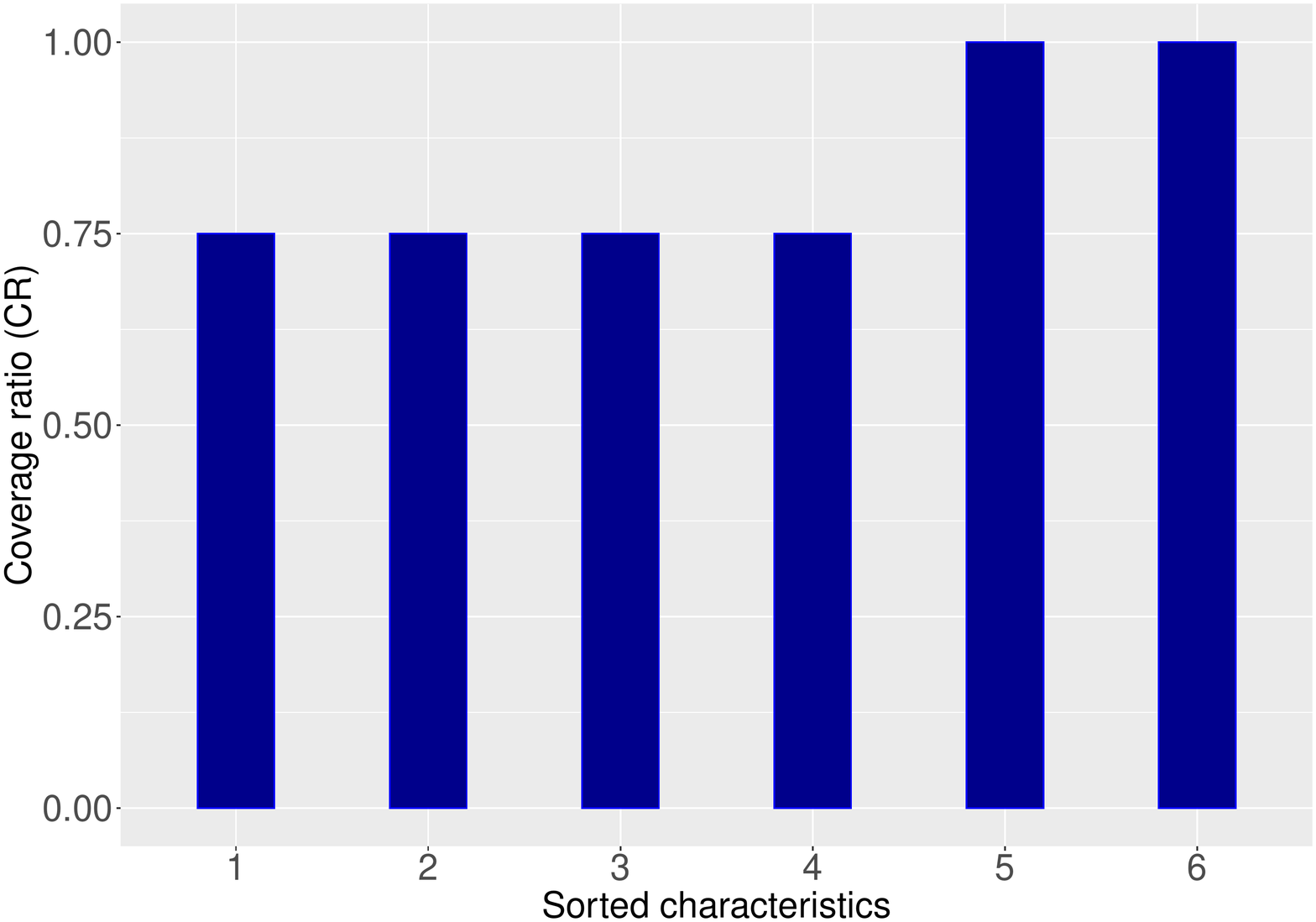}  
\caption{\scriptsize{sorted s2: [0.75, 0.75, 0.75, 0.75, 1, 1]}}
  \label{lexicoC}
\end{subfigure}
\begin{subfigure}{.5\textwidth}
  \centering
  \includegraphics[width=.9\linewidth]{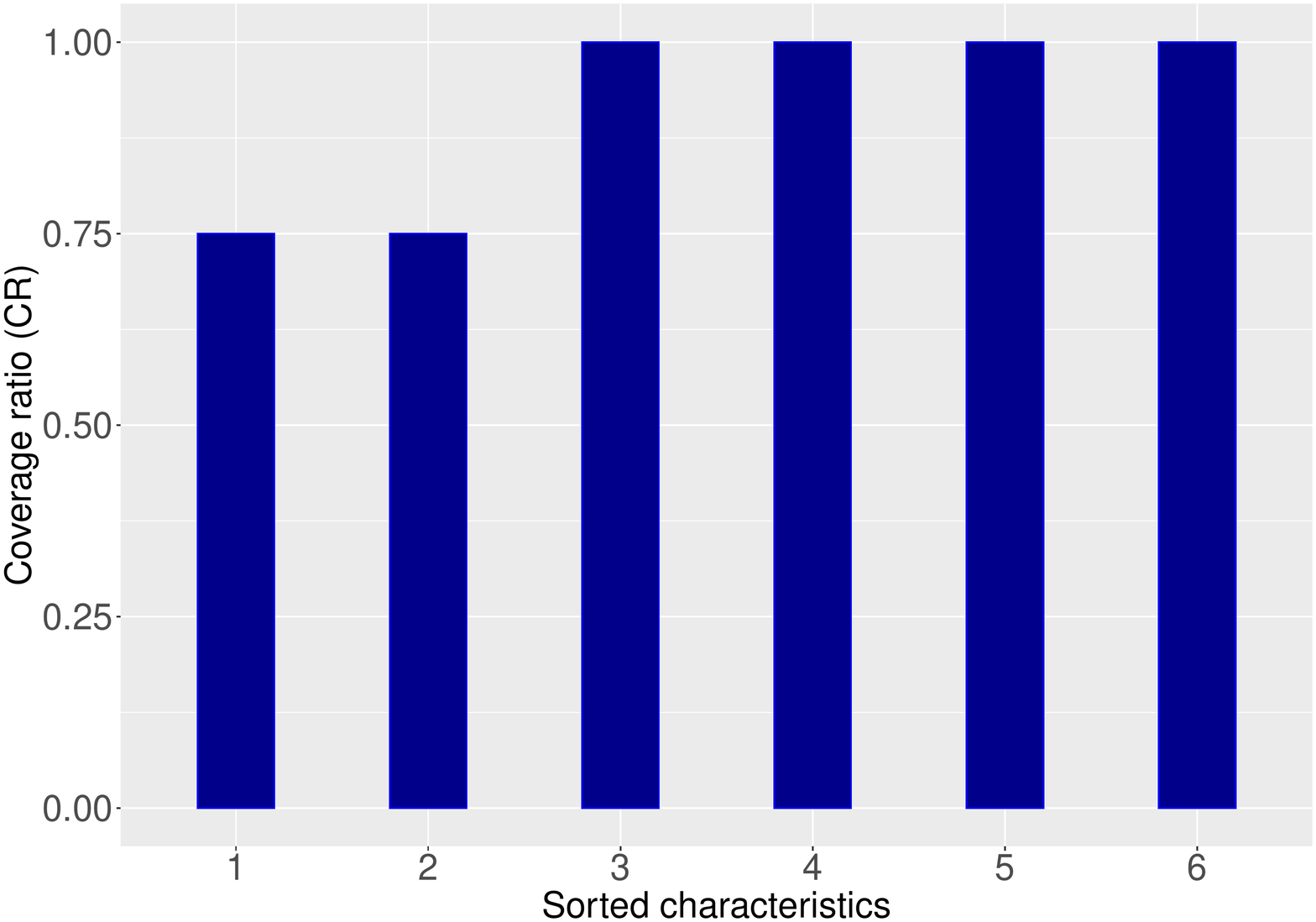}  
\caption{\scriptsize{sorted s3: [0.75, 0.75, 1, 1, 1, 1]}}
  \label{lexicoD}
\end{subfigure}
\caption{Illustrative example for lexicographic maximin approach based on solutions in Fig~\ref{maximin}}
\label{lexico}
\end{figure}

\textbf{Definition 1}. (\textit{Leximin ordering}): Let $x$ and $y$ denote two vectors in $\mathbb{R}^m$. Let $x^{\uparrow}$ and $y^{\uparrow}$ denote these same vectors where each element has been rearranged in a non-decreasing order. According to the leximin ordering:
\begin{itemize}
\item the vector $x$ leximin-dominates $y$ (written $x \succ _{leximin} y$) if and only if $\exists i \in \{1,....,m\}$ such that: $\forall  j \in \{1,...,i-1\}$, $x_{j}^{\uparrow}$ = $y_{j}^{\uparrow}$ and $x_{i}^{\uparrow}$ $>$ $y_{i}^{\uparrow}$,
\item $x$ and $y$ are indifferent (written $x$ $\sim_{leximin}$ $y$ ) if and only if $x^{\uparrow}$ = $y^{\uparrow}$.
\item  $x\succeq _{leximin} y$ is the case where $x \succ _{leximin} y$ or $x$ $\sim_{leximin}$ $y$ .
\end{itemize}

 As for the bi-objective setting, we consider the problem of finding the set of solutions to the SARP that offer the best trade-off between total route duration (we call it duration in the rest of the paper) and coverage ratio maximization. The lexicographic ordering is used to compare solutions on coverage ratio vector maximization aspects. We call this problem the leximin-SARP. As the leximin-SARP is a bi-objective problem, we define Pareto optimality in the leximin-SARP as follows:

\textbf{Definition 2}. \textit{(Dominance and Pareto optimality in the leximin-SARP): }Let  $c_s$ be the duration of solution $s$ and $l^s$ = ($l_{1}^s$ , . . . , $l_{m}^s$ ) the vector of coverage ratios for solution $s$.
\begin{itemize}
\item Let $s$ and $s^\prime$ represent two solutions of the leximin-SARP. A solution $s$ dominates a solution $s^\prime$ iff: $c_s \leq c_{s^\prime}$ and $l^s \succeq _{leximin} l^{s^\prime}$ and either $c_s < c_{s^\prime}$ or $l^s \succ _{leximin} l^{s^\prime}$.
\item $s$ is a Pareto optimal solution iff no other solution dominates $s$.

\end{itemize}

\begin{center}
\begin{figure}[]
\begin{center}
\includegraphics[scale=.50]{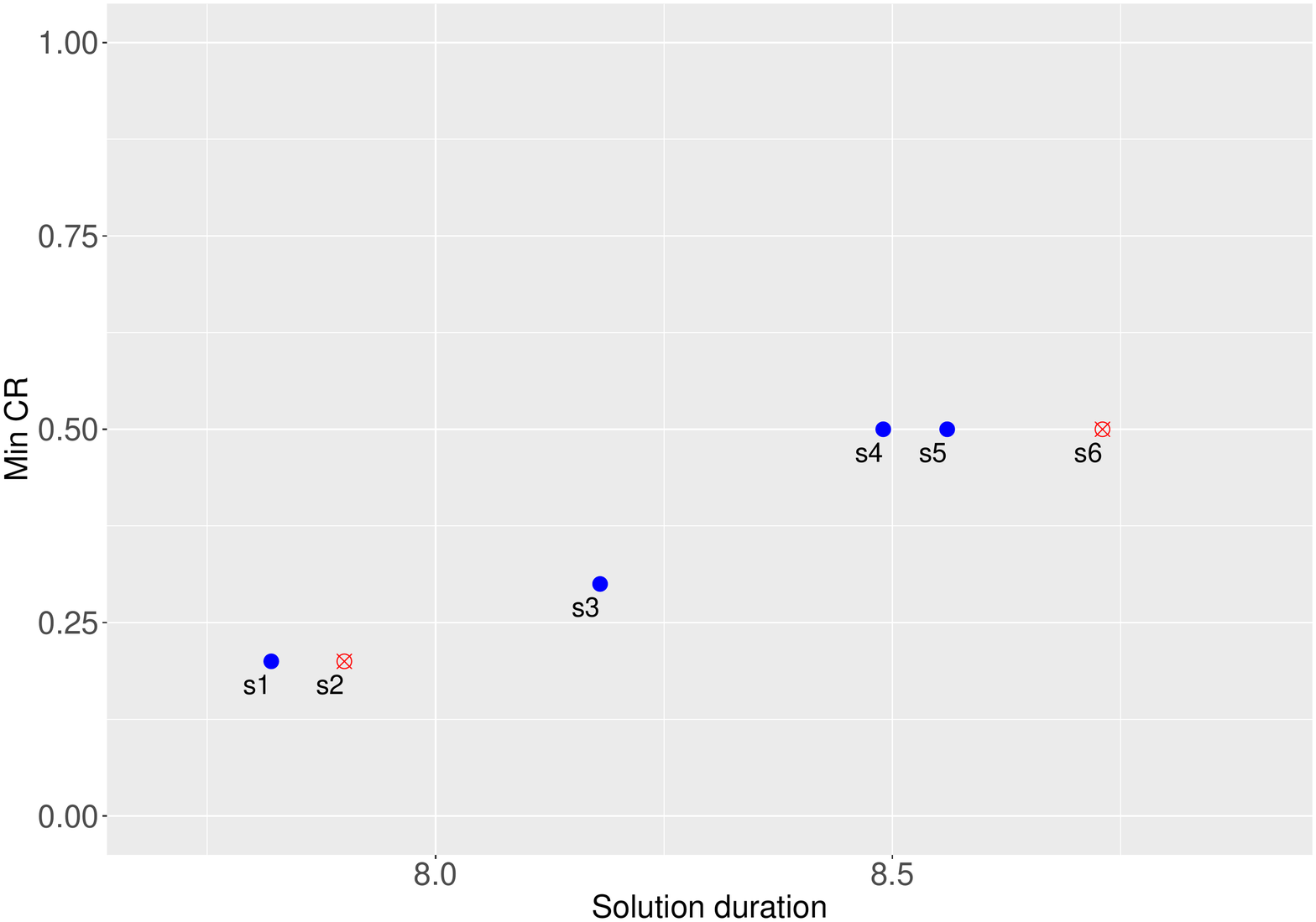}
\end{center}
\caption{Two-dimension visualization of solutions provided in Table~\ref{tableexample} in which solutions s2 and s6 are dominated by other solutions.}
\label{example}
\end{figure}
\end{center}

To illustrate the concepts of dominance and Pareto optimality,  we present the next  example in Table~\ref{tableexample}, including a set of six solutions according to their resulted durations and coverage ratios (CRs). In this example $s_2$ is dominated by $s_1$ since $s_1$ has a lower duration and leximin-dominates $s_2$. Similarly, $s_6$ is dominated by $s_5$. $s_4$, however, is not dominated by $s_5$ since although $s_5$ leximin-dominates $s_4$, $s_4$ has a lower duration. Fig~\ref{example} provides a two-dimension solution visualization for this example. Note that although in lexicographic ordering we always compare two vectors of values instead of two points, in two-dimensional visualization we only plot the max-min value. The information in such a figure should always be interpreted along with solutions' information, i.e., vector of coverage ratios and solution durations (e.g., Table~\ref{tableexample} and Fig~\ref{example}).

\begin{table}[]
  \centering
  \caption{illustrative example for the dominance and Pareto optimality}
    \begin{tabular}{c|r|rrrr}
    \toprule
    \multicolumn{1}{c}{Solution} & \multicolumn{1}{c}{Duration} & \multicolumn{4}{c}{sorted CRs} \\
    \midrule
    $s_1$    & 7.82  & 0.20  & 0.30  & 0.33  & 0.44 \\
    $s_2$    & 7.90  & 0.20  & 0.25  & 0.33  & 0.44 \\
    $s_3$    & 8.18  & 0.30  & 0.30  & 0.33  & 0.44 \\
    $s_4$    & 8.49  & 0.50  & 0.50  & 0.50  & 0.56 \\
    $s_5$    & 8.56  & 0.50  & 0.50  & 0.56  & 0.56 \\
    $s_6$    & 8.73  & 0.50  & 0.50  & 0.54  & 0.57 \\
    \bottomrule
    \end{tabular}%
  \label{tableexample}%
\end{table}%

\section{Solution method}
\label{SolutionMethod}
The SARP,  as a variant of the TOP, is a NP-hard problem \citep{balcik2017site} and consequently the leximin-SARP as well. In order to compute the complete set of Pareto optimal solutions, we develop a bi-objective metaheuristic, that is based on the multi-directional local search (MDLS) framework proposed by \citet{tricoire2012multi}, to approximate the set of Pareto optimal solutions for  the leximin-SARP. We describe the general structure of MDLS  in Section~\ref{mdls}. The components of the operators used in the MDLS are explained in Section~\ref{LNSComponents}.

\subsection{ Multi-directional local search}
\label{mdls}
MDLS is a metaheuristic framework for multi-objective optimization, which generalizes the concept of single-objective local search to multiple objectives.  To perform MDLS with the leximin-SARP as a bi-objective problem, a specific local search should be defined for each objective (duration and leximin). Therefore, we develop a parametric adaptive large neighborhood search algorithm, called $ALNS_i$ , where parameter $i$ represents the objective function being optimized: $ALNS_1$ aims at minimizing durations and  $ALNS_2$ aims at maximizing the leximin value. ALNS is a metaheuristic that iteratively destroys and repairs an incumbent solution and has been successfully applied to various transportation problems \citep{shaw1998using,ropke2006adaptive,pisinger2019large}. We explain the components of the ALNS algorithm in Section~\ref{LNSComponents}. By developing ALNS for both objectives, we can apply MDLS.  

MDLS considers an archive of solutions, $F$, which is preserved and updated via the solution process.  MDLS maintains $F$ as a non-dominated set by performing non-dominated union operations \citep{tricoire2012multi}. The non-dominated union of two sets performs the union of these two sets while removing the elements that are dominated by at least one element of these two sets. Algorithm~\ref{alg1} outlines the MDLS steps.  At each iteration, three steps are performed i) a solution $x$ is randomly selected from $F$, each solution having equal probability, ii) for each objective, single-objective local search is performed on $x$,  and iii) $F$ is updated with solutions obtained by single-objective local search \citep{tricoire2012multi}.  We initialize the archive $F$  with the solution of the cheapest insertion heuristic adapted from \citet{pisinger2007general}.

\begin{algorithm}[]
1: pre-condition: F is a non-dominated set\\
2: \textbf{repeat}\\
3: \hspace{.5cm}$x \gets$ select a solution $(F)$\\
4: \hspace{.5cm}$\textbf{for}\hspace{.1cm} k \in \{1,2\} \hspace{.1cm} \textbf{do}$\\
5: \hspace{.8cm}$F \leftarrow F \cup_{\preceq} ALNS_k(x)$\\
6: \hspace{.5cm}$\textbf{end for}$\\
7: \textbf{until} timeLimit is reached\\
8: \textbf{return} $F$\\
\caption{Multi-directional local search (MDLS) framework proposed  by \citet{tricoire2012multi}}
\label{alg1}
\end{algorithm}

\subsection{ALNS components}
\label{LNSComponents}

 Within the designed MDLS framework, for each objective, we define several \emph{destroy} and \emph{repair} operators. Then, for each objective, at each iteration, we randomly choose a destroy and repair from the set of operators for that objective and generate a new solution.
  As in the leximin-SARP, it is not possible to visit all sites; the sites removed from one solution by the destroy operator are added to the list of unassigned (or not visited) sites. The repair operator then selects and inserts some of the unassigned sites in the partially destroyed solution until no further insertion is possible due to time constraints. Below, we explain the destroy and repair operators for each objective.
  
\subsubsection{Duration objective: destroy and repair operators}
Classical LNS destroy and repair operators for the VRP with the single objective of cost minimization are defined in detail within the literature, e.g., \citet{pisinger2007general}. In this study, for the duration objective, we implement some of these destroy and repair operators. The applied removal operators are random removal, worst removal, and related removal. In the following, we describe each briefly:

\begin{itemize}
    \item Random removal: this operator diversifies the search by randomly selecting sites to be removed from the current solution. We repeat this process until the number of removed sites is equal to the destroy quantity $q$ \citep{pisinger2007general}.
    \item Worst removal: in this operator the $q$ sites with the worst contribution to the
    duration objective value are removed.  The contribution of one site is calculated by the difference in the solution's total route duration with and without this site.  This operator randomizes the selected sites to avoid situations where the same sites are removed repetitively. In particular, we sort the sites by descending duration then we choose a random number $y$ in the interval [0,1). The site to be removed is $L[y^{p_{worst}}|L|]$ ; where $L$ represents the array of sorted sites and  the parameter $p_{worst}$ controls the degree of randomization  \citep{pisinger2007general}.  
   \item Related removal: this operator selects a random site $i$ and removes it as well as its $q-1$  most related sites, aiming to shuffle these similar sites around and create new, perhaps better, solutions. We define relatedness $R(i, j)$ of two sites $i$ and $j$ by: $R(i, j)$ = $t_{ij}$, where $t_{ij}$ is the routing cost (based on travel distance). The randomization works similar to the worst removal operator, with parameter $p_{related}$ controlling the degree of randomization \citep{pisinger2007general}. 

\end{itemize}
The  two repair operators used for the duration objective work as follows:
\begin{itemize}
    \item Cheapest insertion heuristic: this operator calculates the minimum insertion cost for every unassigned site.  The insertion cost is determined by subtracting the current solution's total duration without the inserted site from the duration of the solution with the inserted site. Afterward, we insert the site with the minimum cost difference in the current solution. We continue this process until no further insertion is possible, i.e., maximum allowed time ($T_{max}$) is reached \citep{pisinger2007general}.
    \item k- regret heuristic: this heuristic aims to improve the cheapest insertion heuristic by incorporating a  look-ahead information when selecting the site to insert. The site with the biggest regret value is inserted, i.e., the site that has the highest difference between the cost of the insertion into the best tour and the insertions into the k best tours \citep{pisinger2007general}. A k- regret heuristic calculates the part with the greatest cost difference between the cheapest and the k-1 next cheapest insertions.  The selected site is inserted at its minimum cost position. We consider  2- regret heuristics and 3-regret heuristics as two separate repair operators.
\end{itemize}

\subsubsection{Leximin objective: destroy and repair operators}

For the leximin objective, we use the following two destroy  operators:
\begin{itemize}
\item \textit{Random removal}:  this operator is the same as the one used for the duration objective and aims to diversify the search by randomly selecting sites to be removed from the current solution.  This process is repeated until the number of removed sites is equal to the destroy quantity $q$.
\item \textit{Worst min removal}: in this operator the $q$ sites with the worst contribution to the minimum coverage ratio of the solution are removed. That is, this operator starts removing sites that decrease the least the minimum coverage ratio of the current solution. The randomization works similar to the duration objective, controlled by parameter $p_{worst}$.
\end{itemize}
Two repair operators have been designed to guide the search towards leximin-optimal solutions:
\begin{itemize}
\item \textit{Highest max-min  insertion}: this operator prioritizes the selection of sites with the highest minimum coverage ratio of characteristics to be inserted first. That is, we first insert (if feasible) the site with the highest contribution to the max-min value of the current solution. We insert each unassigned site within a feasible position that minimizes the total travel time. In order to discriminate between insertions with a same contribution to the max-min objective, we provide two versions of this operator:

\begin{itemize}
\item \textit{Highest random max-min:} in case there are equivalent max-min values while sorting the sites, the candidate site for insertion is randomly selected among those sites with best max-min value. We continue this process until no further insertion is possible. 
\item \textit{Highest max-min with priority to the duration}: in this method, among the sites with equivalent max-min values,  we select the one that yields the lowest increase to the duration objective function. We continue this process until no further insertion is possible.
   
\end{itemize} 

\item \textit{Highest Leximin  insertion}: in this operator, we prioritize the selection of sites with the highest contribution to the leximin objective. More precisely,  we compare the respective vectors of coverage ratios resulting from inserting the unassigned sites. This is in contrast to the \textit{Highest max-min  insertion} which is a proxy to the leximin objective and considers only the max-min value. Therefore, this operator is computationally expensive. We evaluate the performance of this operator in Section~\ref{QualityOfConfigs}. 

\end{itemize} 

We follow the adaptive weight adjustment described in \citep{ropke2006adaptive} to automatically adjust the weights of different operators using statistics from earlier iterations. According to this method, we keep track of the success rate of each operator, which measures how well the operator has performed recently. Within the scope of our study, an operator is successful when using this operator leads to an update in the non-dominated set of solutions. We divide the entire search into a number of segments, each corresponding to the same number of MDLS iterations. We calculate new weights using the recorded scores after each segment using this formula: 
$\omega _{i,j+1}=\omega _{i,j}(1-r)+r\frac{\pi_i}{\theta_i}$, 
where $\omega _{i,j}$ represents the weight of operator $i$ used in segment $j$,  $\pi_i$ is the score of operator $i$ obtained during the last segment, and $\theta_i$ is the number of times we have
attempted to use operator $i$ during the last segment. The parameter $r$ (reaction factor) controls how fast the weight adjustment algorithm reacts to changes in the effectiveness of the operators. This is similar to the scheme used by~\cite{ropke2006adaptive}.

\section{Computational results}
\label{Computational}
The designed algorithm is implemented in C++. All computational experiments are carried out on a 64-bit Windows Server with two 2.3 GHz Intel Core CPUs and 12 GB RAM. 
We describe our test instances and three different configurations in Section~\ref{InstanceDescription}. Parameter settings are presented in Section~\ref{ParameterSettings}. The quality of the Pareto front approximation is explained in Section~\ref{QualityOfConfigs}. We analyse the equivalent solutions for the max-min objective function in Section~\ref{EquivalentSolutionsSection} and compare our results with the max-min values reported in\citet{balcik2017site} in  Section~\ref{AnalysisOnThemax-min}.

\subsection{Instance and configuration description}
\label{InstanceDescription}

We evaluate our designed MDLS framework on the \citet{balcik2017site} instances with 25, 50, 75, and 100 nodes. These instances are generated by modifying Solomon’s 100-node random (R) and random-clustered (RC) instances \citep{balcik2017site}. 
Within these instances, there are 12  characteristics. The maximum allowed duration for each team ($T_{max}$) is set between two and eight, and the number of assessment teams ($K$)  ranges between two and six. As a result, there are 48 instances.  Table~\ref{instances} shows the characteristics of the SARP instances for our experiments. 

We run MDLS on each instance 10 times with two different run times, using three different configurations called \textit{all}, \textit{leximin} and \textit{max-min}. The difference between these configurations is whether we use the costly operator of \textit{Highest Leximin insertion} or not. In fact, the three configurations include all removal operators as well as all duration-oriented operators. The \textit{all} configuration includes both Highest max-min insertion and \textit{Highest Leximin insertion}. The \textit{leximin} configuration includes only \textit{Highest Leximin  insertion} and the \textit{max-min} configuration includes only \textit{Highest max-min  insertion}.

\begin{table}[]
  \caption{Characteristics of the SARP instances for our experiments. The columns shows the number of sites in each instance, the type of network, number of teams, the maximum allowed duration for each route and run times.}
  \label{instances}
  \centering
  \scalebox{0.8}{
\begin{tabular}{cccccc}
\hline
\multicolumn{1}{l}{instance} & \multicolumn{1}{l}{N-type/K/T_{max}} & \multicolumn{1}{l}{time limit (sec.)} & \multicolumn{1}{l}{instance} & \multicolumn{1}{l}{N-type/K/T_{max}} & \multicolumn{1}{l}{time limit (sec.)} \\ \hline
R1                             & 25\_R/2/2                             & \multirow{6}{*}{30, 90}             & RC1                           & 25\_RC/2/2  & \multirow{6}{*}{30, 90}             \\
R2                             & 25\_R/2/3                             &                                    & RC2                           & 25\_RC/2/3 &                                    \\
R3                             & 25\_R/2/4                             &                                    & RC3                           & 25\_RC/2/4 &                                    \\
R4                             & 25\_R/3/2                              &                                    & RC4                           & 25\_RC/3/2 &                                    \\
R5                             & 25\_R/3/3                              &                                    & RC5                           & 25\_RC/3/3 &                                    \\
R6                             & 25\_R/3/4                              &                                    & RC6                           & 25\_RC/3/4                             &                                  \\\hline
R7                             & 50\_R/3/3                              & \multirow{6}{*}{60, 180}            & RC7                           & 50\_RC/3/3  & \multirow{6}{*}{60, 180}            \\
R8                             & 50\_R/3/4                             &                                    & RC8                           & 50\_RC/3/4   &                                    \\
R9                             & 50\_R/3/5                             &                                    & RC9                           & 50\_RC/3/5  &                                    \\
R10                            & 50\_R/4/3                             &                                    & RC10                          & 50\_RC/4/3   &                                    \\
R11                            & 50\_R/4/4                            &                                    & RC11                          & 50\_RC/4/4                               &                                    \\
R12                            & 50\_R/4/5                            &                                    & RC12                          & 50\_RC/4/5 &                                    \\\hline
R13                            & 75\_R/3/3                            & \multirow{6}{*}{120, 360}           & RC13                          & 75\_RC/3/3 & \multirow{6}{*}{120, 360}           \\
R14                            & 75\_R/3/4                              &                                    & RC14                          & 75\_RC/3/4  &                                    \\
R15                            & 75\_R/3/6                              &                                    & RC15                          & 75\_RC/3/6 &                                    \\
R16                            & 75\_R/5/3                              &                                    & RC16                          & 75\_RC/5/3 &                                    \\
R17                            & 75\_R/5/4                               &                                    & RC17                          & 75\_RC/5/4 &                                    \\
R18                            & 75\_R/5/6                               &                                    & RC18                          & 75\_RC/5/6 &                                    \\\hline
R19                            & 100\_R/3/4                              & \multirow{6}{*}{240, 720}           & RC19                          & 100\_RC/3/4 & \multirow{6}{*}{240, 720}           \\
R20                            & 100\_R/3/6                             &                                    & RC20                          & 100\_RC/3/6 &                                    \\
R21                            & 100\_R/3/8                            &                                    & RC21                          & 100\_RC/3/8 &                                    \\
R22                            &100\_R/6/4                            &                                    & RC22                          & 100\_RC/6/4  &                                    \\
R23                            & 100\_R/6/6                            &                                    & RC23                          & 100\_RC/6/6 &                                    \\
R24                            & 100\_R/6/8                            &                                    & RC24                          & 100\_RC/6/8 &  \\ \hline                                
\end{tabular}}
\end{table}
\subsection{Parameter settings}
\label{ParameterSettings}
 We performed preliminary computations to test our MDLS performance under different parameter settings, which led to the following parameters: The ruin quantity $q$  used in the destroy operators is randomly selected within bounds between at least one site and at most 30 percent of the number of sites in the current solution. The other parameters are set as follows: ($p_{related}$, $p_{worst}$, $r$) = (5, 3, 0.1). For the purpose of adaptive weight adjustment, each segment lasts 100 MDLS iterations.

\subsection{ Quality of the Pareto front approximation}
\label{QualityOfConfigs}

To analyze the impact of using different configurations on the quality of the Pareto front approximation, we follow the approach described in \citet{lehuede2020lexicographic}. Accordingly,  we provide a reference set for each instance. The reference set for each instance is the non-dominated union of the sets returned by the algorithm over 60 runs, including all three MDLS configurations, two different time limits, ten runs per setting. We compare the quality of each set by looking at the percentage of found solutions from the reference set. In addition, we analyze the percentage of reference solutions found within a  1\%, 2\%, and 3\% distance of a solution from the assessed set. We call solution $s_1$  within a $\alpha$\% distance of another solution $s_2$ if, when the duration and all coverage ratios of $s_2$ are multiplied by 1 + $\alpha/100$  and 1 - $\alpha/100$ respectively,  then $s_1$ dominates this transformed solution. The boxplots provided in Appendix~\ref{boxplots} (Figures~\ref{R25} to~\ref{RC100}) represent the percentage of reference solutions found and those within a  1\%, 2\%, and 3\% distance of a solution from the reference set.  The variability in the boxplots is derived from the results returned by 10 runs for each MDLS configuration,  each instance, and each run time.

We see that the three tested configurations of MDLS are close to the reference set for almost all instances.  Moreover, we observe that all three configurations provide relatively similar approximations of the Pareto set in most cases. Therefore, the simple construction heuristics of the \textit{max-min} configuration are generally enough to reach an acceptable set. Note that, although \textit{all}  and \textit{leximin} configurations include  the computationally expensive operator of \textit{Highest Leximin insertion}, they do not result in better solution sets in most cases.
 
\begin{figure}[h]
 \centering
  \includegraphics[width=1\linewidth]{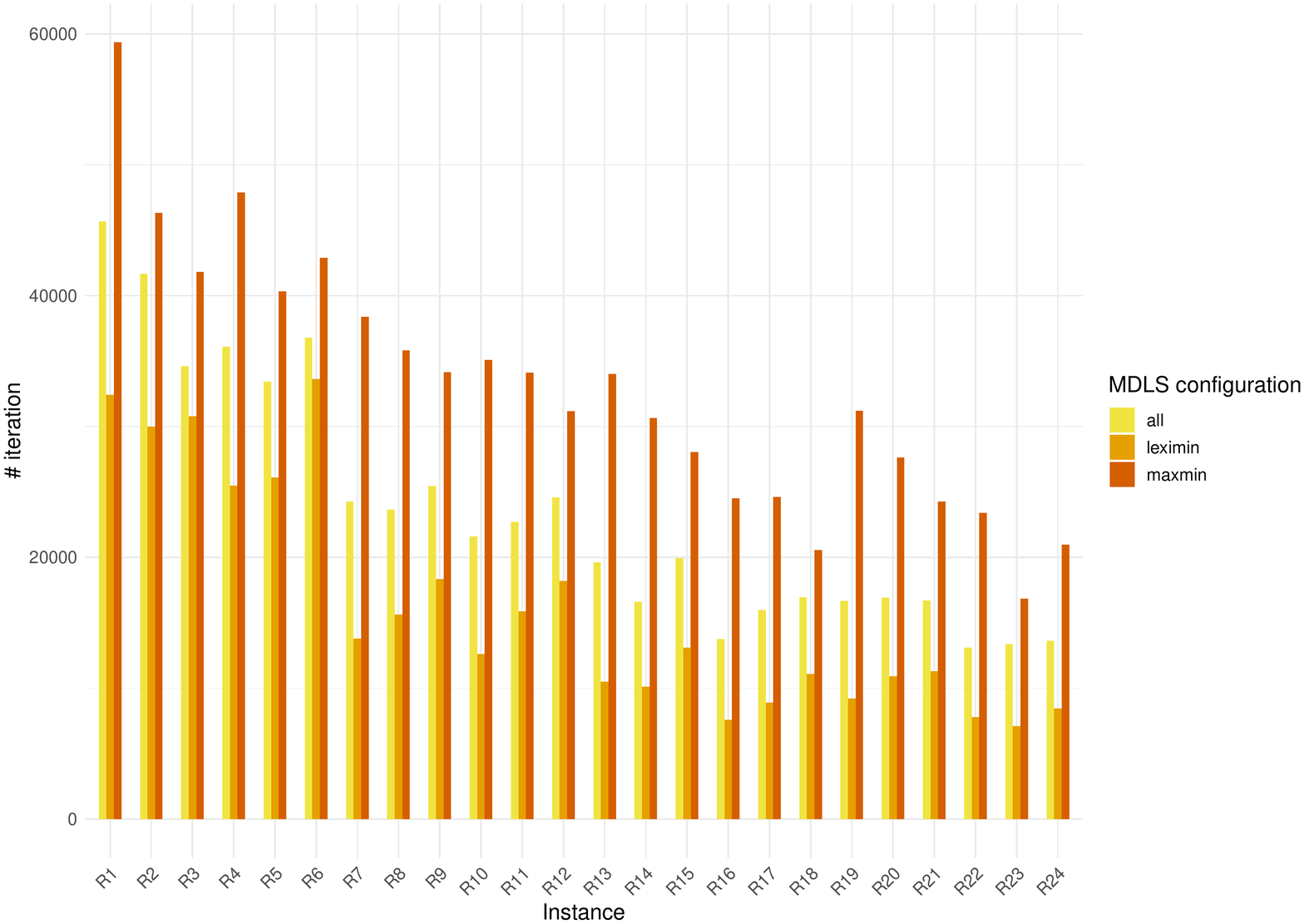}  
 \caption{Average number of MDLS iterations for each configuration on each instance.}
    \label{Riteration}
   \end{figure}

\begin{figure}[h]
  \centering
   \includegraphics[width=1\linewidth]{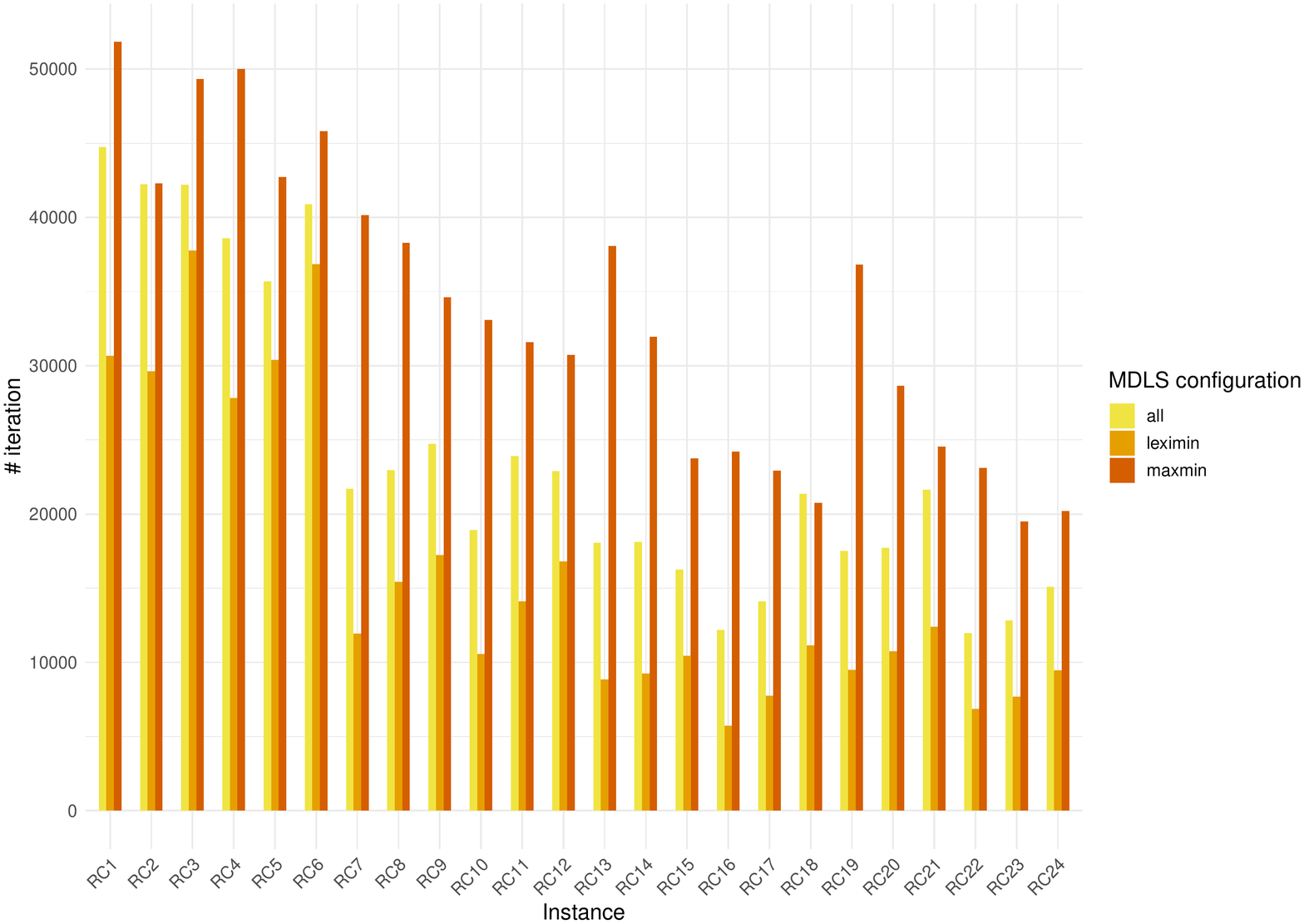}  
   \caption{Average number of MDLS iterations for each configuration on each instance.}
     \label{RCiteration}
   \end{figure}

 To evaluate this impact, we compare the average number of iterations performed by the MDLS algorithm with the various configurations for each instance. Figure~\ref{Riteration} and Figure~\ref{RCiteration} show the average number of iterations performed for each instance. As foreseen, integrating the \textit{Highest Leximin insertion} operator significantly impacts the number of performed iterations.  By looking at both boxplots in Appendix~\ref{boxplots} and Figures~\ref{Riteration} and \ref{RCiteration}, we see that for some instances, more iterations do not necessarily result in better solutions. For example, in instance RC19, even though the \textit{max-min} configuration performs almost three times as many iterations as in  \textit{leximin}, using the \textit{Highest Leximin insertion} heuristic is important to obtain better solutions.


 \subsection{Equivalent solutions for the max-min objective function} 
 \label{EquivalentSolutionsSection}
 The motivation for studying the leximin-SARP is that solutions with a same max-min value are not necessarily equivalent, and that it is worth exploring them. 
 In fact, the max-min approach fails to distinguish between solutions with the same minimum values (coverage ratio in this study). 
 Worse, optimising max-min first and duration second in a lexicographic fashion, as is done in~\cite{balcik2017site}, may bias the search towards solutions that are good with the max-min objective but comparatively bad with the leximin objective, since removing visits reduces the secondary objective.
 We call solutions that have the same max-min value \emph{max-min-equivalent}. For example, Table~\ref{RC7MAximinEquivalent} shows 12 solutions with the same max-min value of 0.571 for the instance RC8. These solutions, however, are different in terms of leximin objective. Solution $s5$ in Table~\ref{RC7MAximinEquivalent} is the one that is reported in \citet{balcik2017site} as the best max-min value. Table~\ref{RC7MAximinEquivalent} lists solutions that provide a better vector of coverage ratios than $s5$ while incurring a modest increase in duration. These solutions include more sites in the assessment plan, which is beneficial since each additional visit increases the coverage ratio associated with at least one characteristic.  In particular, solution $s12$ and $s1$ provide the same max-min values but considerably different solutions in terms of leximin objective. We also note that our bi-objective approach also allows us to find solutions that would dominate $s5$ if using the (max-min, duration) lexicographic objective: $s1, s2, s3, s4$. 
 
\begin{table}[]
\centering
\caption{max-min-equivalent solutions for instance RC8}
 \scalebox{0.8}{
\begin{tabular}{c|c|cccccccccccc}
solution & duration & \multicolumn{12}{c}{sorted coverage   ratios}                                                 \\ \hline 
s1       & 10.4 & 0.571 & 0.571 & 0.571 & 0.583 & 0.600 & 0.600 & 0.600 & 0.607 & 0.654 & 0.667 & 0.727 & 0.750 \\  
s2       & 10.5 & 0.571 & 0.571 & 0.600 & 0.636 & 0.643 & 0.643 & 0.667 & 0.692 & 0.700 & 0.720 & 0.722 & 0.750 \\
s3       & 10.6 & 0.571 & 0.600 & 0.625 & 0.636 & 0.643 & 0.643 & 0.654 & 0.667 & 0.679 & 0.700 & 0.720 & 0.750 \\
s4       & 10.7 & 0.571 & 0.636 & 0.643 & 0.643 & 0.643 & 0.650 & 0.654 & 0.667 & 0.667 & 0.700 & 0.720 & 0.750 \\
\textbf{s5}       & \textbf{10.8} & \textbf{0.571} & \textbf{0.643} & \textbf{0.643} & \textbf{0.650} & \textbf{0.700} & \textbf{0.714} & \textbf{0.722} & \textbf{0.727} & \textbf{0.731} & \textbf{0.750} & \textbf{0.750} & \textbf{0.760} \\
s6     & 10.9 & 0.571 & 0.643 & 0.650 & 0.700 & 0.714 & 0.727 & 0.731 & 0.750 & 0.750 & 0.750 & 0.778 & 0.800 \\
s7       & 11.1 & 0.571 & 0.643 & 0.667 & 0.679 & 0.692 & 0.700 & 0.714 & 0.727 & 0.750 & 0.750 & 0.760 & 0.800 \\
s8       & 11.2 & 0.571 & 0.692 & 0.700 & 0.714 & 0.714 & 0.714 & 0.722 & 0.727 & 0.750 & 0.750 & 0.800 & 0.800 \\
s9       & 11.3 & 0.571 & 0.700 & 0.714 & 0.714 & 0.727 & 0.731 & 0.750 & 0.750 & 0.750 & 0.778 & 0.800 & 0.840 \\
s10      & 11.6 & 0.571 & 0.700 & 0.714 & 0.714 & 0.778 & 0.786 & 0.800 & 0.808 & 0.840 & 0.875 & 0.909 & 0.917 \\
s11      & 11.7 & 0.571 & 0.714 & 0.727 & 0.731 & 0.750 & 0.750 & 0.750 & 0.778 & 0.786 & 0.786 & 0.800 & 0.880 \\
s12      & 11.9 & 0.571 & 0.727 & 0.750 & 0.750 & 0.750 & 0.769 & 0.786 & 0.786 & 0.786 & 0.800 & 0.833 & 0.920
\end{tabular}}
\label{RC7MAximinEquivalent}
\end{table} 

\begin{figure}[]
 \centering
  \includegraphics[width=.9\linewidth]{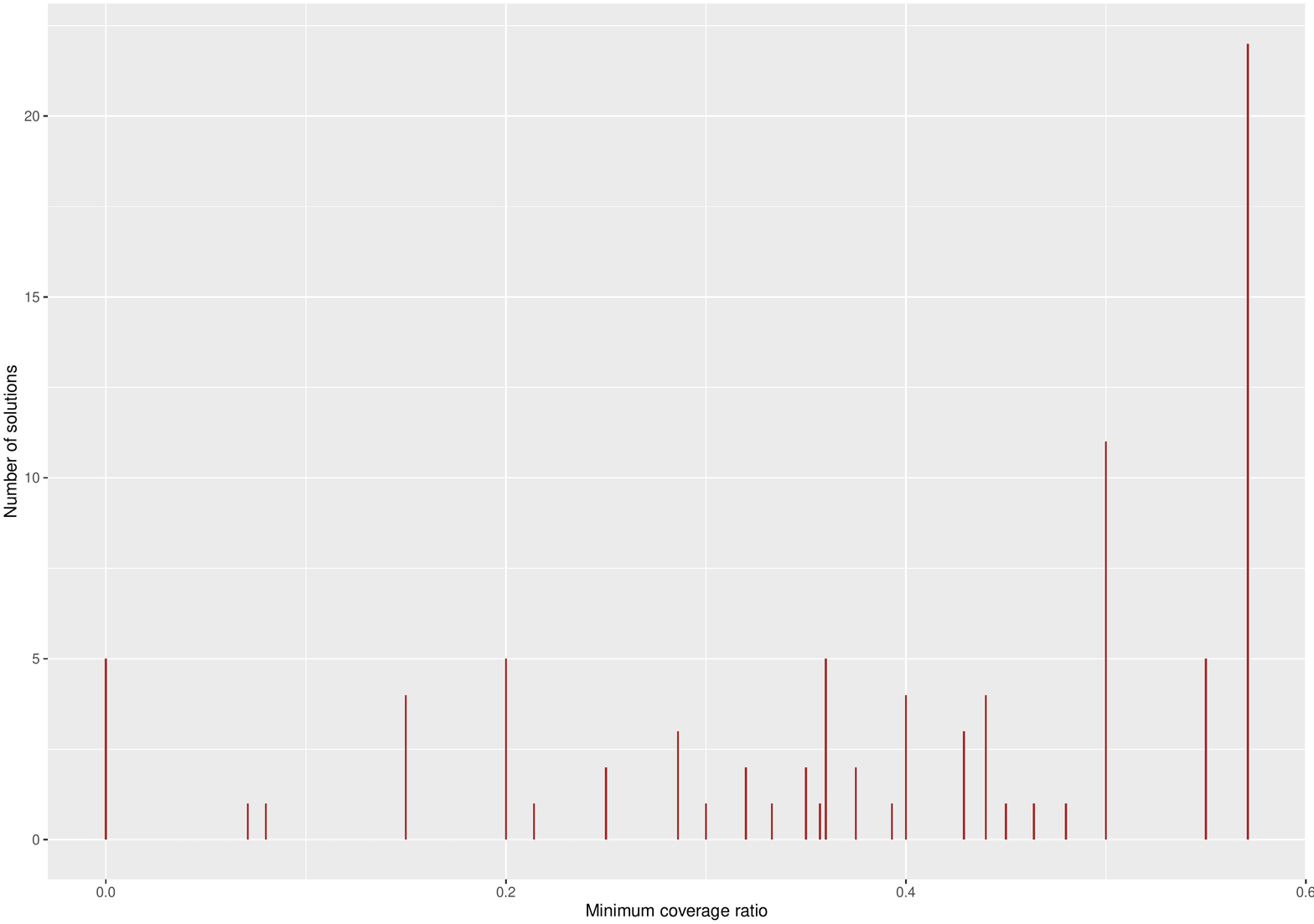}  
  \caption{Distribution analysis for minimum coverage ratio of non-dominated solutions  in the reference set for instance RC8}
  \label{RC8histogram}
  \end{figure}

Figure~\ref{RC8histogram} depicts the histogram of the distribution of the number of max-min-equivalent solutions in the reference set for the instance RC8.  We observe that in this instance, we have 22  distinct solutions with the same minimum coverage ratio of 0.571. We provide similar histograms for other instances in  Appendix~\ref{DistributionAnalysis}. We can see  that for example, instance RC16 has more than 70 distinct solutions with the  same minimum coverage ratio of 0.125.  We also observe that this trend is not contained to highest max-min values; e.g., the instance R11 has more than 9 max-min-equivalent solutions for different values of minimum coverage ratio. Note that these histograms include only non-dominated solutions in the solution sets found for the leximin-SARP. 

 \subsection{Comparison with max-min values reported in previous approaches} 
 \label{AnalysisOnThemax-min}
 
\begin{figure}[H]
 \centering

  \includegraphics[width=1\linewidth]{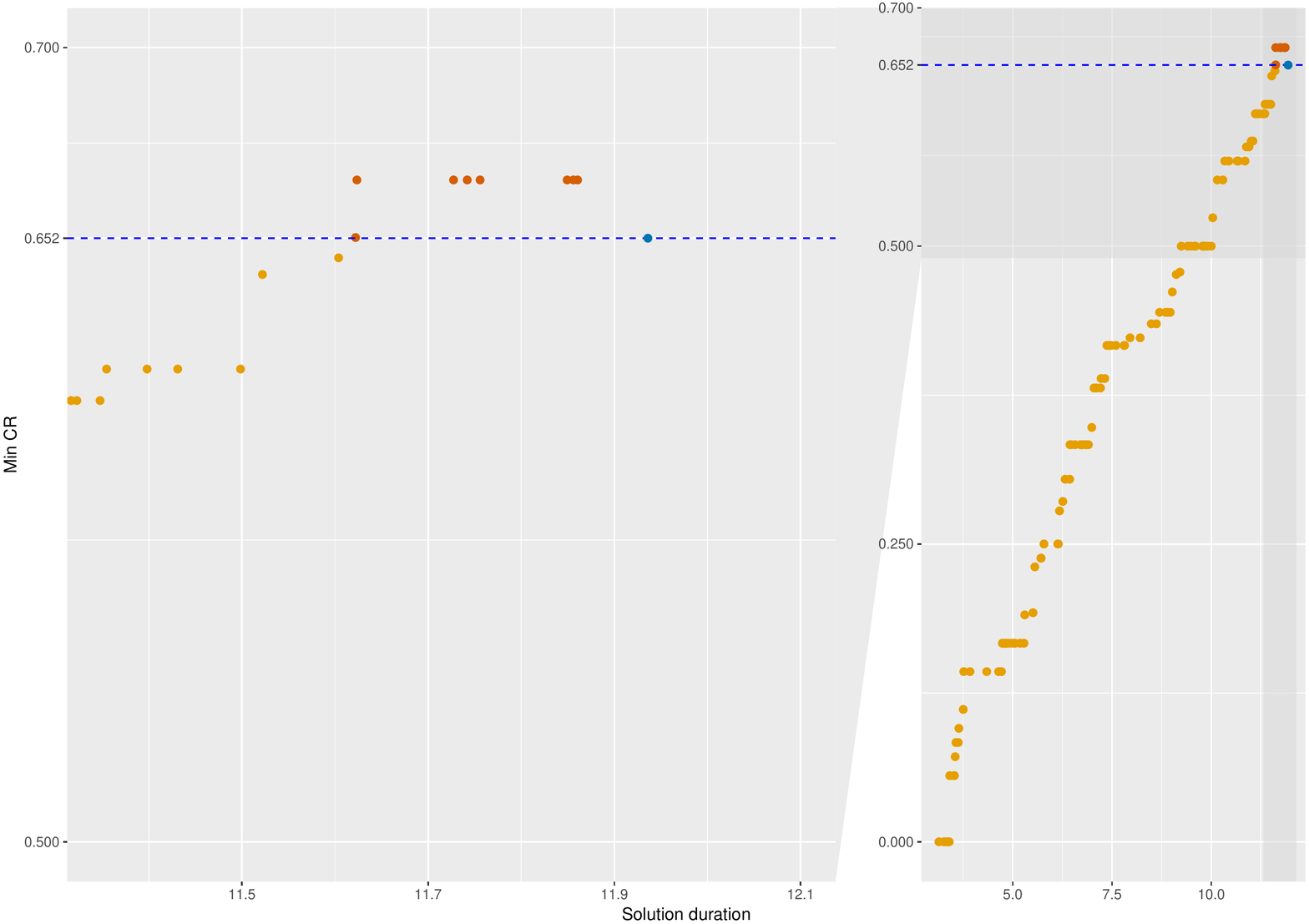}  
  \caption{Union of Non-dominated solution sets of different MDLS variants for the instance R8 (left part is the zoomed-in area)}
  \label{R8solutionsDominateBalcik2017}
  \end{figure}

\begin{figure}[H]
 \centering
   \includegraphics[width=1\linewidth]{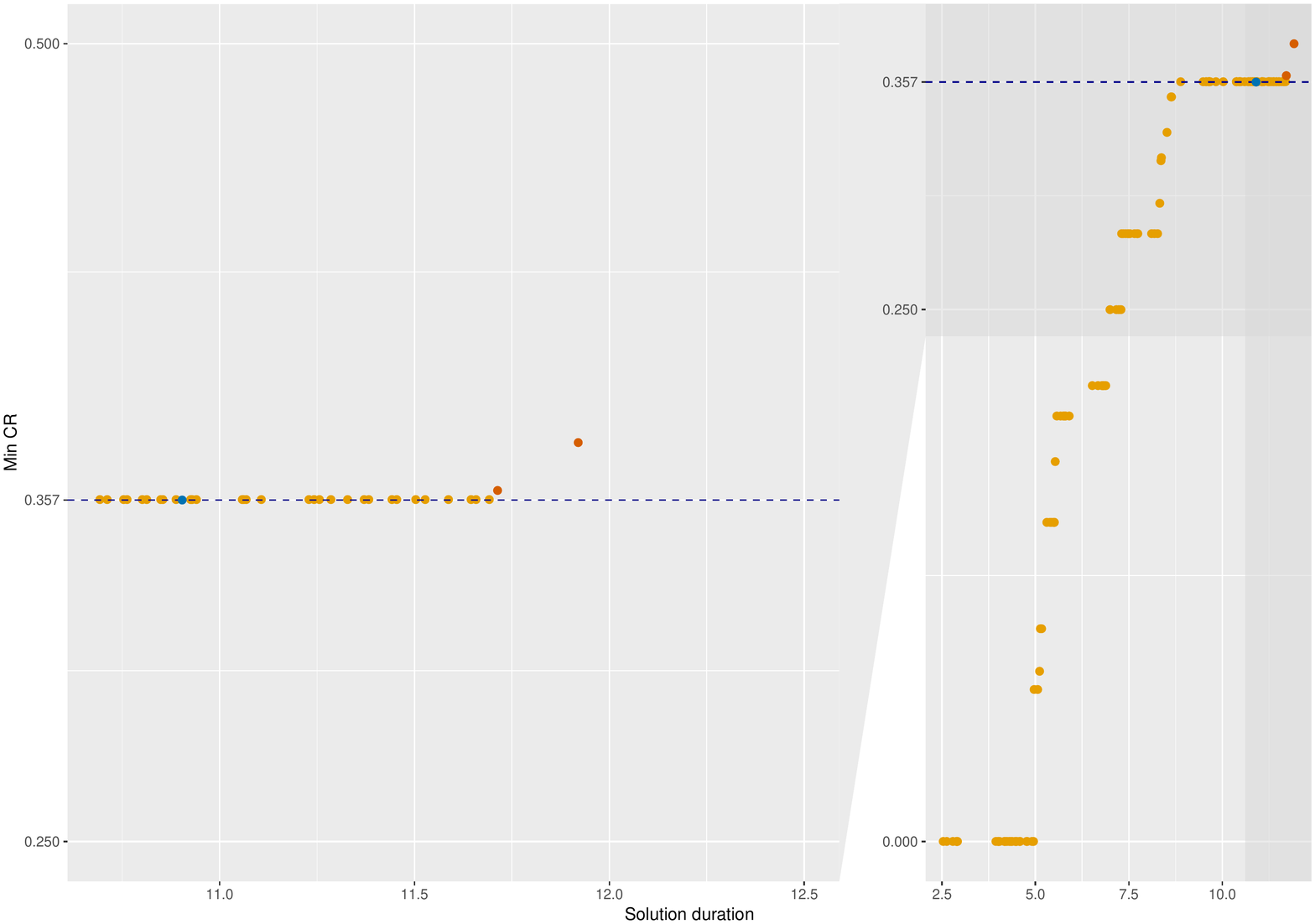}  
  \caption{Union of Non-dominated solution sets of different MDLS variants for the instance RC10 (left part is the zoomed-in area)}
  \label{R21solutionsDominateBalcik2017}
  \end{figure}

The primary objective function in  \citet{balcik2017site}  is maximizing the minimum coverage ratio. Minimizing total route duration is considered as a secondary objective in a lexicographic fashion. Our results show that this method can lead to solutions that are dominated by the leximin approach. Figure~\ref{R8solutionsDominateBalcik2017} provides a graphical representation of the non-dominated set of solutions obtained by the MDLS algorithm for the instance R8. The solution duration is represented in the horizontal axis, and the minimum coverage ratio is represented in the vertical axis. In this figure, the blue point represents the max-min solution reported in \citet{balcik2017site} which is dominated by other solutions obtained by MDLS. Moreover, we find new solutions that have a higher values in terms of max-min than the previously best known value. These points are represented in red (above the dashed line representing the previously best known max-min value). Note that we do not have access to the full vector of coverage ratios of the solutions reported in \citet{balcik2017site}. Hence, the comparisons are based on the reported minimum value of coverage ratios and their corresponding durations. For example,  Figure~\ref{R21solutionsDominateBalcik2017} shows the case for which we cannot say that the solutions obtained by the MDLS algorithm dominate the one in \citet{balcik2017site}. However, we know that there are two solutions that leximin-dominate it (the two solutions with a higher max-min coverage ratio).
The graphical representations of other instances in which we obtained such results are provided in Appendix~\ref{Union of Non-dominated}.

We now look at how many distinct solutions we find with the same max-min value as the one reported in \citet{balcik2017site}, counting all solutions from reference sets.  Table~\ref{numberofhighermax-min} provides that information for each instance separately.

It is important to mention that  in most of the runs the designed MDLS  finds the same max-min coverage ratios and, in some cases, an even higher value, than the previously best known value. In particular, for long run times, \textit{max-min}, \textit{all}, and \textit{leximin} configurations find a solution with a max-min value at least as good as the previously best known one in 87, 85 and 68 percent of the runs respectively (see Appendix~\ref{percentageOfBestmax-mins} for the detailed percentages for each instance).  This is in addition to providing a set of trade-off solutions. We also note that the bi-objective MDLS algorithm uses comparable CPUs (two 2.0 GHz Intel Xeon versus two 2.3 GHz Intel Core). MDLS produces approximation sets which includes the previously best-known max-min values, in most of the cases within shorter run times, for all but small sizes.

\begin{table}[]
\caption{Number of minimum coverage ratio in the reference set that are same or higher than those reported in the \citet{balcik2017site}}
 \scalebox{0.7}{
   \begin{tabular}{|c|ccc|}
    \toprule
    \multicolumn{1}{|c|}{\multirow{2}[2]{*}{Instance}} & \multicolumn{1}{c}{\multirow{2}[2]{*}{max-min value in Balcik (2017)}} & \multicolumn{1}{c}{\multirow{2}[2]{*}{\# of solutions with same max-min }} & \multicolumn{1}{c|}{\multirow{2}[2]{*}{\# of solutions with higher max-min }} \\
          &       &       &  \\
    \midrule
    R1    & 0.167 & 3     & 0 \\
    R2    & 0.333 & 2     & 0 \\
    R3    & 0.571 & 1     & 0 \\
    R4    & 0.167 & 10    & 0 \\
    R5    & 0.5   & 8     & 0 \\
    R6    & 0.889 & 1     & 0 \\
    \midrule
    R7    & 0.417 & 8     & 0 \\
    R8    & 0.652 & 1     & 7 \\
    R9    & 0.833 & 4     & 0 \\
    R10   & 0.571 & 2     & 1 \\
    R11   & 0.833 & 8     & 0 \\
    R12   & 1.000     & 1     & 0 \\
    \midrule
    R13   & 0.429 & 3     & 0 \\
    R14   & 0.588 & 7     & 3 \\
    R15   & 0.85  & 3     & 2 \\
    R16   & 0.65  & 1     & 1 \\
    R17   & 0.85  & 1     & 3 \\
    R18   & 1.000     & 1     & 0 \\
    \midrule
    R19   & 0.571 & 1     & 0 \\
    R20   & 0.821 & 2     & 5 \\
    R21   & 0.957 & 5     & 5 \\
    R22   & 0.912 & 1     & 3 \\
    R23   & 1.000     & 1     & 0 \\
    R24   & 1.000     & 1     & 0 \\
    \midrule
    RC1   & 0.154 & 1     & 0 \\
    RC2   & 0.438 & 1     & 0 \\
    RC3   & 0.667 & 1     & 0 \\
    RC4   & 0.167 & 3     & 0 \\
    RC5   & 0.667 & 10    & 0 \\
    RC6   & 0.889 & 1     & 0 \\
    \midrule
    RC7   & 0.357 & 8     & 0 \\
    RC8   & 0.571 & 22    & 0 \\
    RC9   & 0.857 & 8     & 0 \\
    RC10  & 0.357 & 50    & 2 \\
    RC11  & 0.846 & 1     & 0 \\
    RC12  & 1.000     & 1     & 0 \\
    \midrule
    RC13  & 0.435 & 1     & 0 \\
    RC14  & 0.435 & 3     & 2 \\
    RC15  & 0.905 & 2     & 2 \\
    RC16  & 0.125 & 76    & 0 \\
    RC17  & 0.75  & 20    & 0 \\
    RC18  & 1.000     & 1     & 0 \\
    \midrule
    RC19  & 0.444 & 6     & 0 \\
    RC20  & 0.774 & 5     & 6 \\
    RC21  & 0.96  & 1     & 2 \\
    RC22  & 0.84  & 5     & 0 \\
    RC23  & 1.000     & 1     & 0 \\
    RC24  & 1.000     & 1     & 0 \\
    \bottomrule
    \end{tabular}}%
  \label{numberofhighermax-min}%
\end{table}

 
\section{Conclusion} 
\label{conclusion}
In this study, we investigate how the lexicographic maximin approach can further improve the coverage ratio concerns in the rapid needs assessment routing problem. We consider the bi-objective optimization of total assessment time and coverage ratio vector maximization and provide an approximation of the Pareto set for leximin–SARP using multi-directional local search. We develop classical neighborhood operators for the duration and coverage ratio vector maximization objectives. These operators are used within the multi-directional local search framework. We show that implementing them improves the quality of the produced solution sets.

Computational experiments also support our premise that different solutions with the same max-min values are not equivalent. We see that for most of the instances, we find several max-min-equivalent solutions that are, in fact, different solutions in terms of the lexicographic maximin approach - including those that offer a better vector of coverage ratios with a modest increase in duration. As a by-product of our study, we see that the designed MDLS, in some cases,  finds solutions with higher max-min values than previously known ones or solutions that dominate the previously best known solutions with regard to the new leximin objective.

Introducing the lexicographic maximin approach to the rapid needs assessment stage is a first step to show that there is an alternative to the classic max-min approach that can express equity in a better way. Future research can focus on using this approach in other contexts of humanitarian logistics, such as last mile distribution problems where fairness to beneficiaries is of vital importance. 
Future research can also consider uncertainty with respect to travel time.  Due to the high level of uncertainty in the condition of the transportation network after the occurrence of a disaster, complete information on travel times may not be available during the rapid needs assessment stage. This uncertainty can be captured using robust optimization or stochastic programming.


\section*{Acknowledgments}

The authors want to thank \citet{balcik2017site} for providing us with the data set used in her study.

\begin{appendices}

\section{Boxplots }

The boxplots provided in this section represent the proportion of solutions in the reference set found within a 0\%, 1\%, 2\%, 3\% distance for each MDLS configuration within short and long time limits.

We observe that the three tested configurations of MDLS do not find more than  15\% of the reference set except for small instances with 25 nodes (see Figures~\ref{R25} to \ref{RC100}-part a). However, for almost all instances and methods, the produced set is indeed close to the reference set. In most of the cases more than 70\% of the returned sets are within 1\% distance from the reference set (see Figures~\ref{R25} to \ref{RC100}-part b). This trend clearly increases by considering 2\% and 3\% distance from the reference set (Figures~\ref{R25} to \ref{RC100}-part c and d). 

These boxplots also show that all three configurations provide relatively similar approximations of the Pareto set in most of the cases. Therefore the results returned by the \textit{max-min} configuration, which is the simplest among the three configurations, are generally enough to reach an acceptable set.
 Nevertheless, in some cases, the \textit{Highest Leximin insertion} heuristic, used in the \textit{leximin} and \textit{all} configurations, is helpful to some degree to find better solutions. For example for instance RC19 (Figure~\ref{RC100}b-d), configurations \textit{leximin} and \textit{all} perform better. In \textit{Highest Leximin insertion} a  vector of coverage ratios is compared, while in \textit{Highest max-min  insertion}  only the max-min value is evaluated. Therefore,  the \textit{max-min} configuration is faster and performs more iterations within the same time budget.

During our computational experiments, we observed that in general, when the average number of node insertions per iteration of the \textit{leximin} objective is very low ($\approx$ less than 4), the \textit{leximin} configuration performs better than the \textit{max-min} configuration. The reason might be that in this situation, there is a limited possibility to insert nodes. Therefore,  when a new node is inserted while using the \textit{max-min} configuration and this insertion is not improving the leximin objective,  there is limited possibility to make up for it later.

\begin{landscape}
\label{boxplots}
\begin{figure}
 \centering
  \includegraphics[width=1\linewidth]{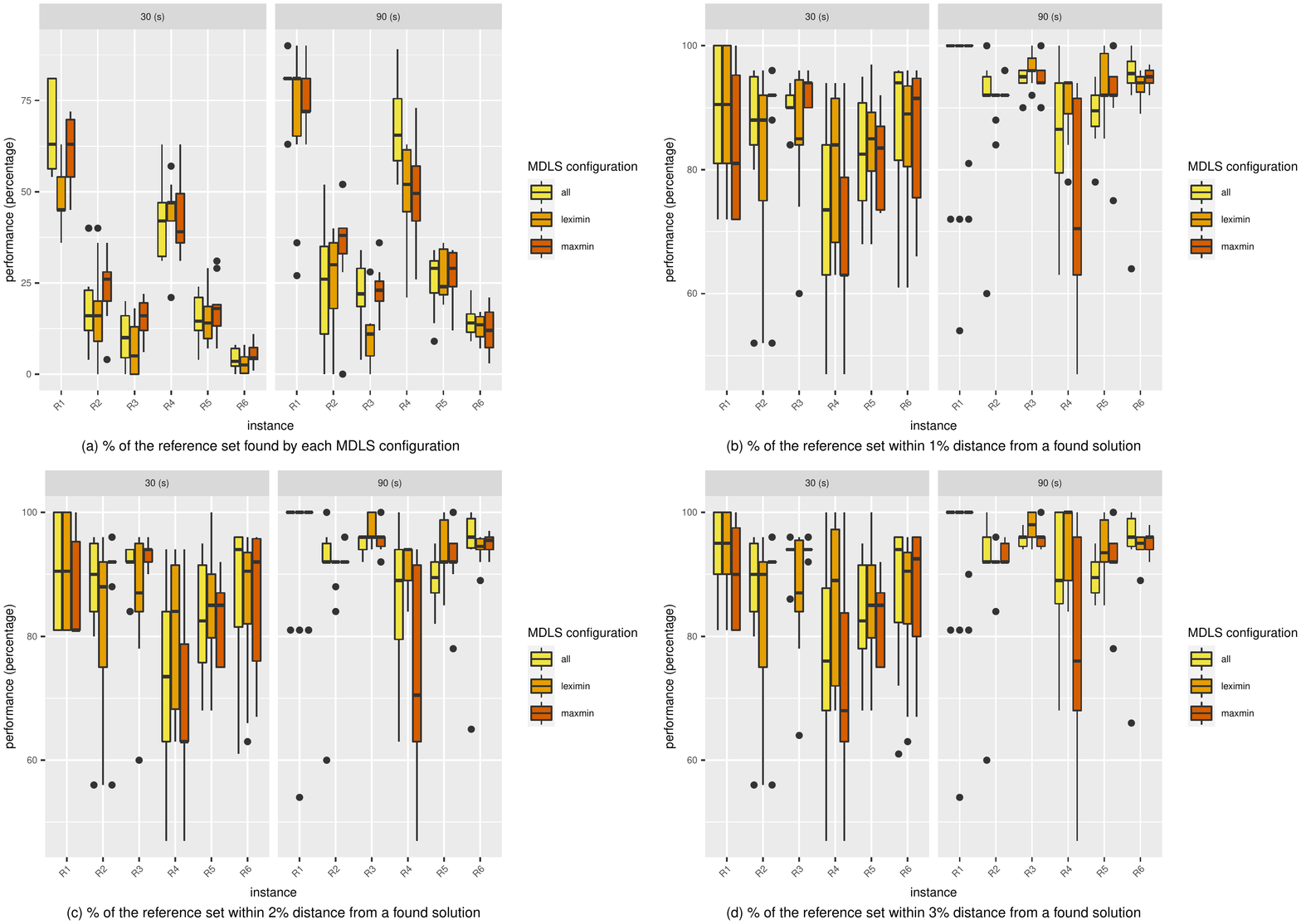}  
  \caption{Proportion of solutions in the reference set found within a 0\%, 1\%, 2\%, 3\% distance (number of sites: 25\_R, run time: 30 and 90 seconds).}
  \label{R25}
  \end{figure}
  
  \begin{figure}
 \centering
  \includegraphics[width=1\linewidth]{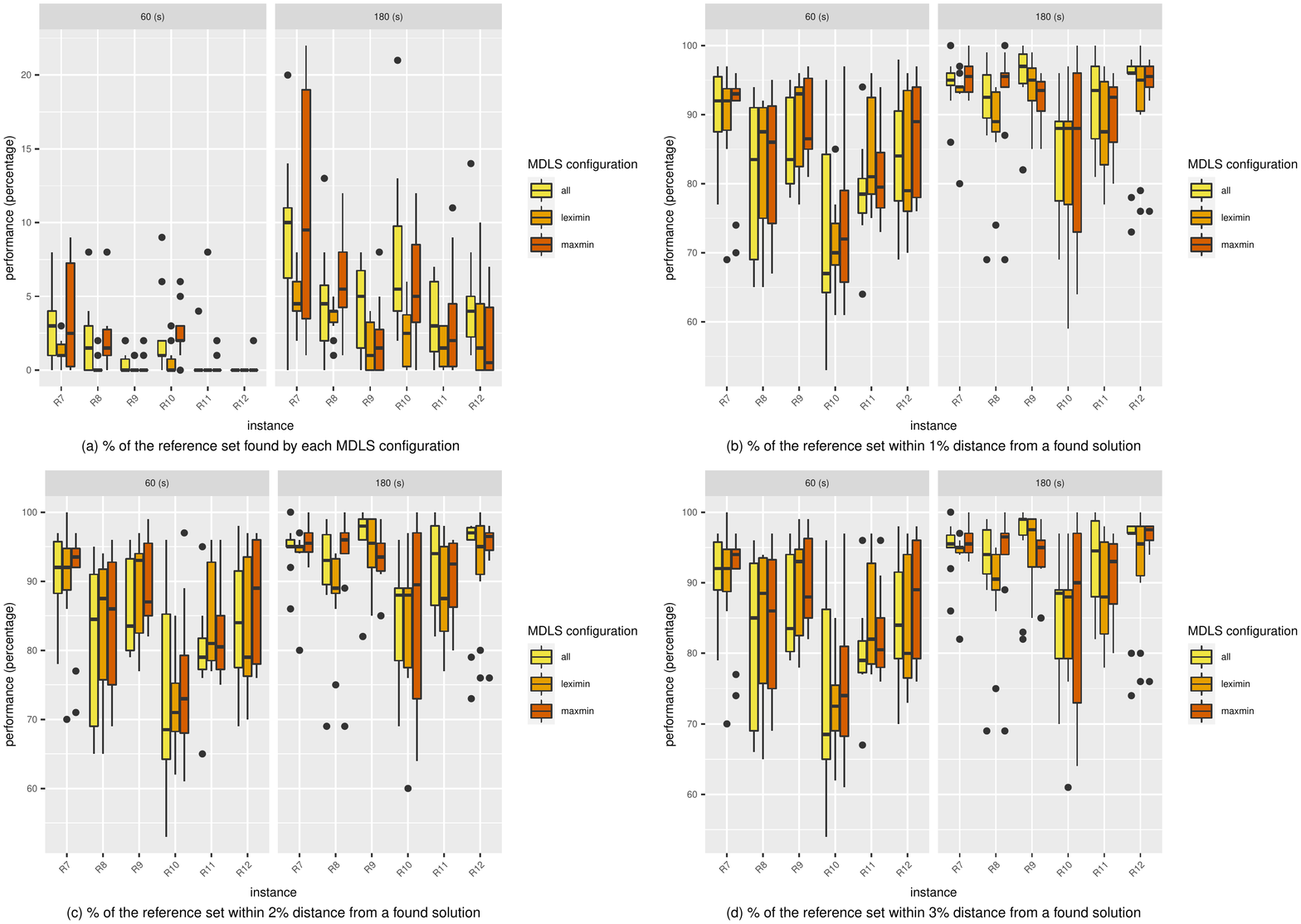}   
  \caption{Proportion of solutions in the reference set found within a 0\%, 1\%, 2\%, 3\% distance (number of sites: 50\_R, run time: 60 and 180 seconds).}
    \label{R50}
  \end{figure}

  \begin{figure}
 \centering
 \includegraphics[width=1\linewidth]{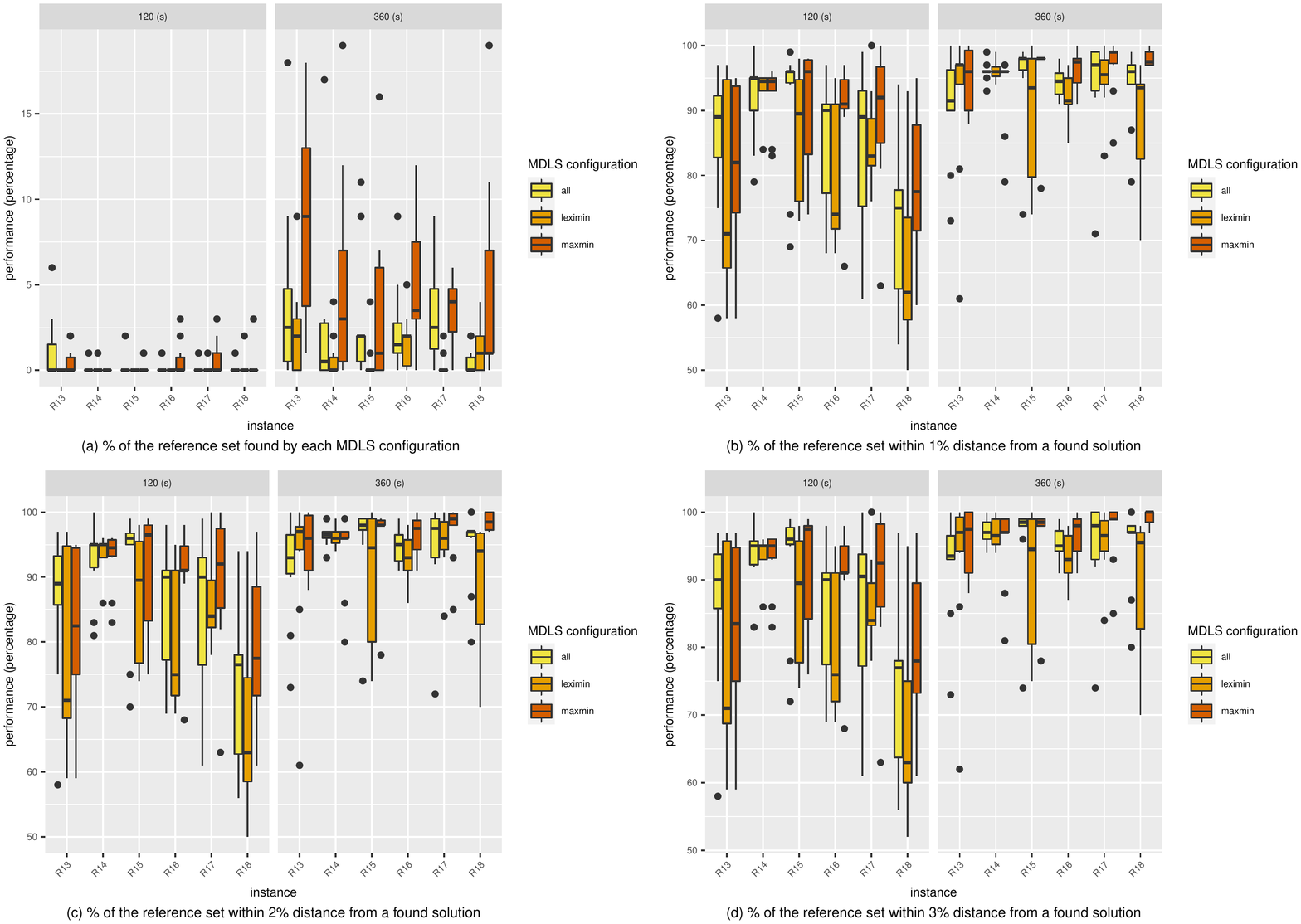}   
  \caption{Proportion of solutions in the reference set found within a 0\%, 1\%, 2\%, 3\% distance (number of sites: 75\_R, run time: 120 and 360 seconds).}
    \label{R75}
  \end{figure}

  \begin{figure}
 \centering
  \includegraphics[width=1\linewidth]{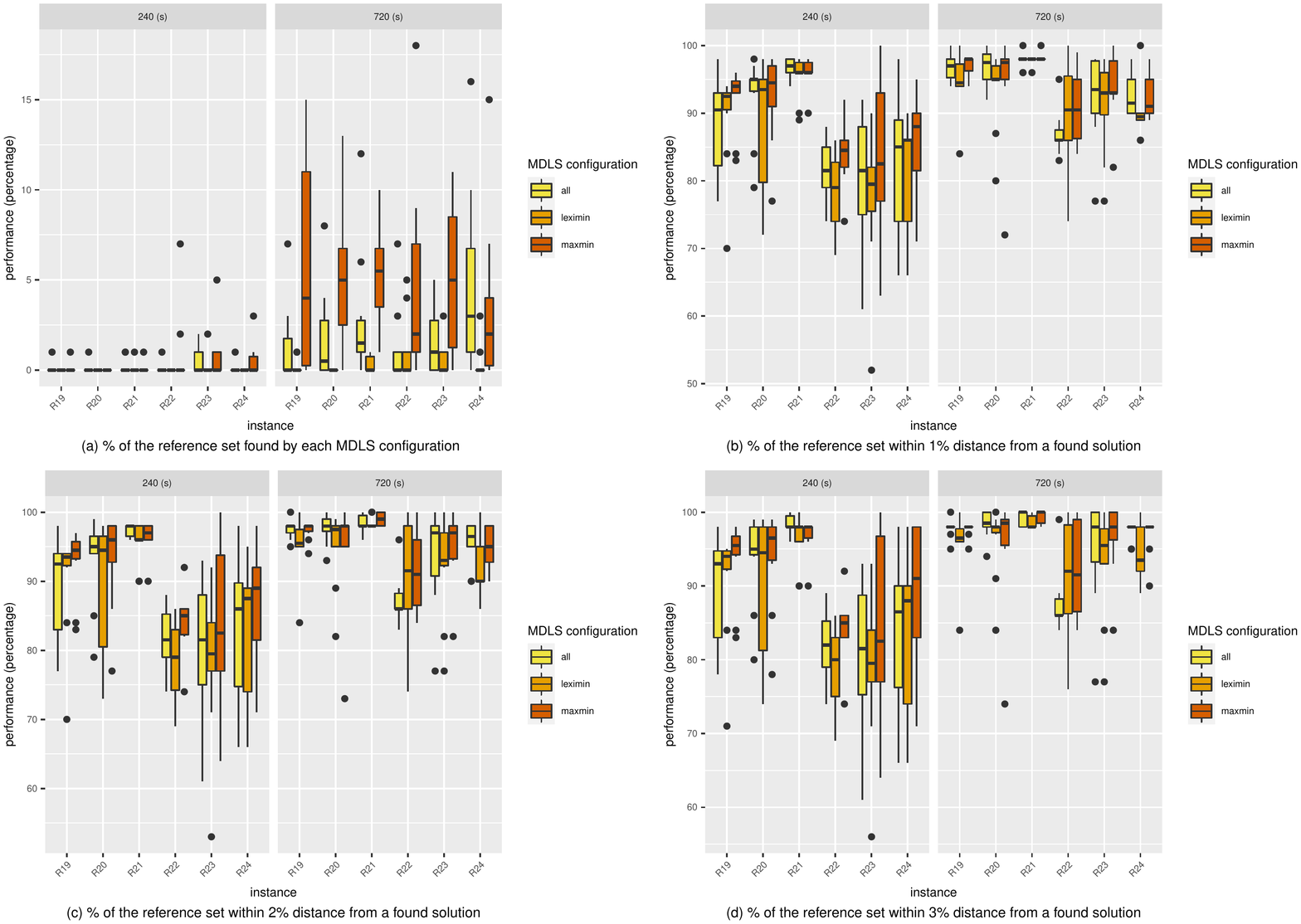}  
  \caption{Proportion of solutions in the reference set found within a 0\%, 1\%, 2\%, 3\% distance (number of sites: 100\_R, run time: 240 and 720 seconds).}
    \label{R100}
  \end{figure}
  
 \begin{figure}
 \centering
  \includegraphics[width=1\linewidth]{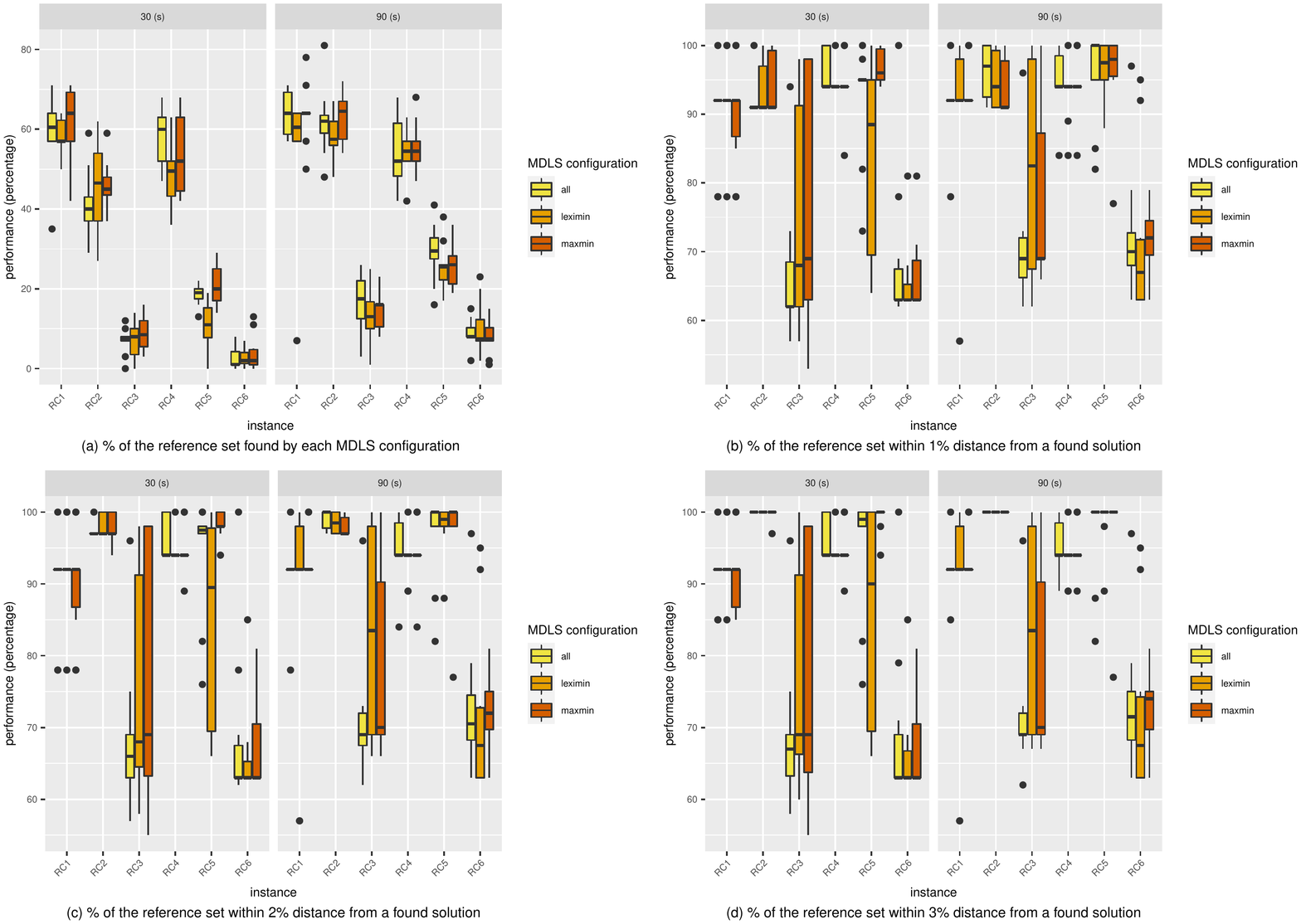}  
  \caption{Proportion of solutions in the reference set found within a 0\%, 1\%, 2\%, 3\% distance (number of sites: 25\_RC, run time: 30 and 90 seconds).}
   \label{RC25}
  \end{figure}
  
  \begin{figure}
 \centering
  \includegraphics[width=1\linewidth]{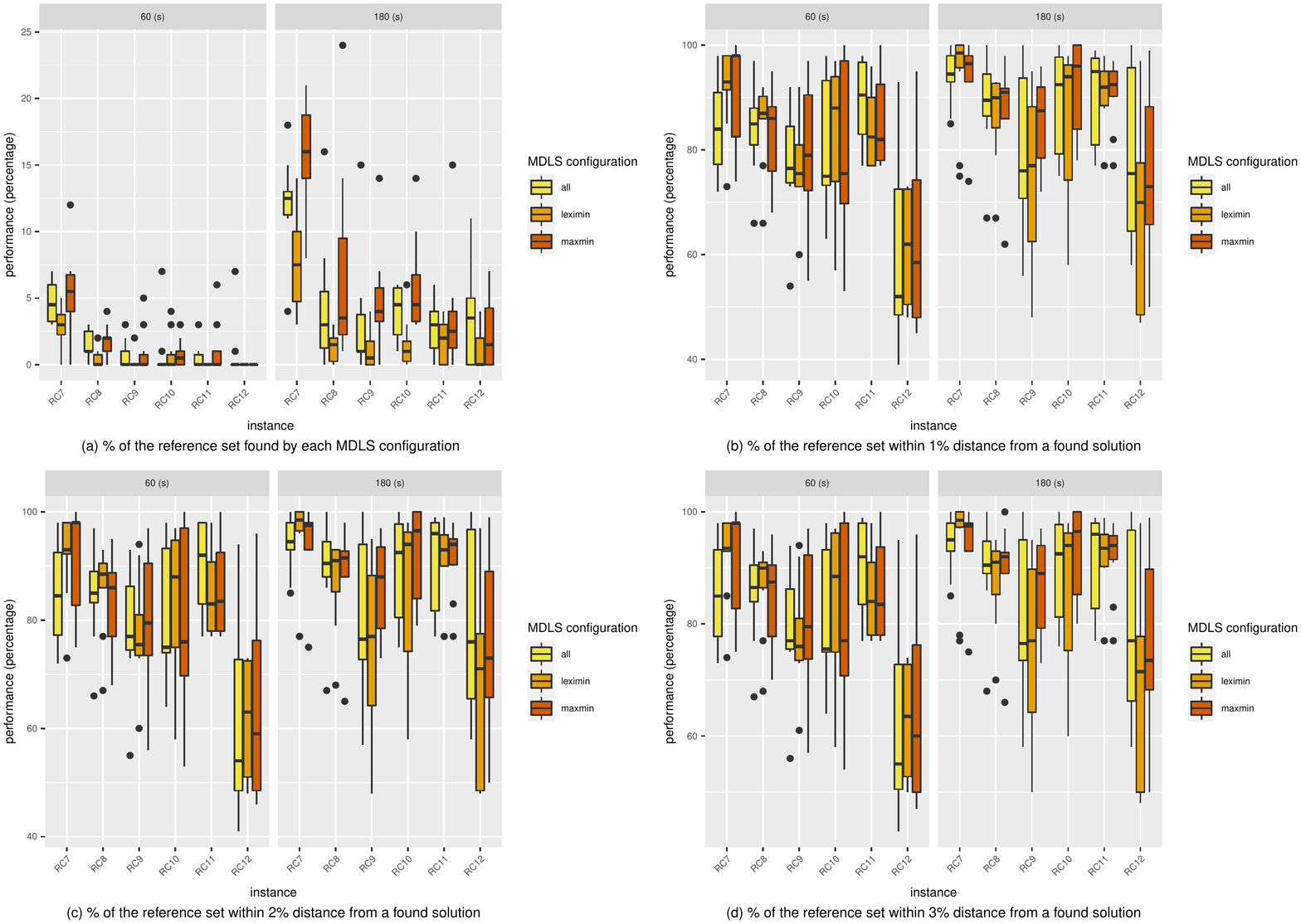}   
  \caption{Proportion of solutions in the reference set found within a 0\%, 1\%, 2\%, 3\% distance (number of sites: 50\_RC, run time: 60 and 180 seconds).}
    \label{RC50}
  \end{figure}

  \begin{figure}
 \centering
 \includegraphics[width=1\linewidth]{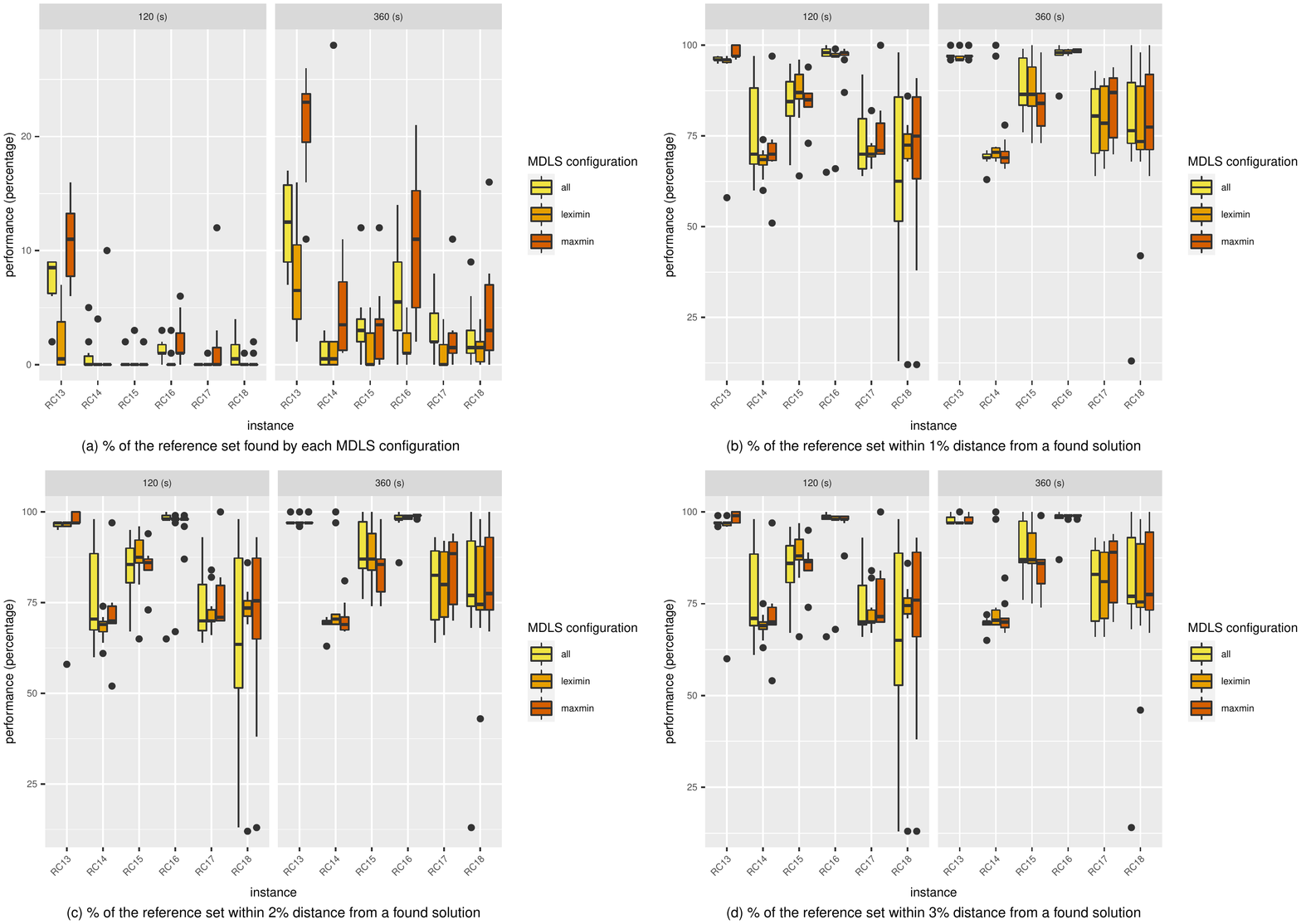}   
  \caption{Proportion of solutions in the reference set found within a 0\%, 1\%, 2\%, 3\% distance (number of sites: 75\_RC, run time: 120 and 360 seconds).}
    \label{RC75}
  \end{figure}
  
   \begin{figure}
 \centering
  \includegraphics[width=.95\linewidth]{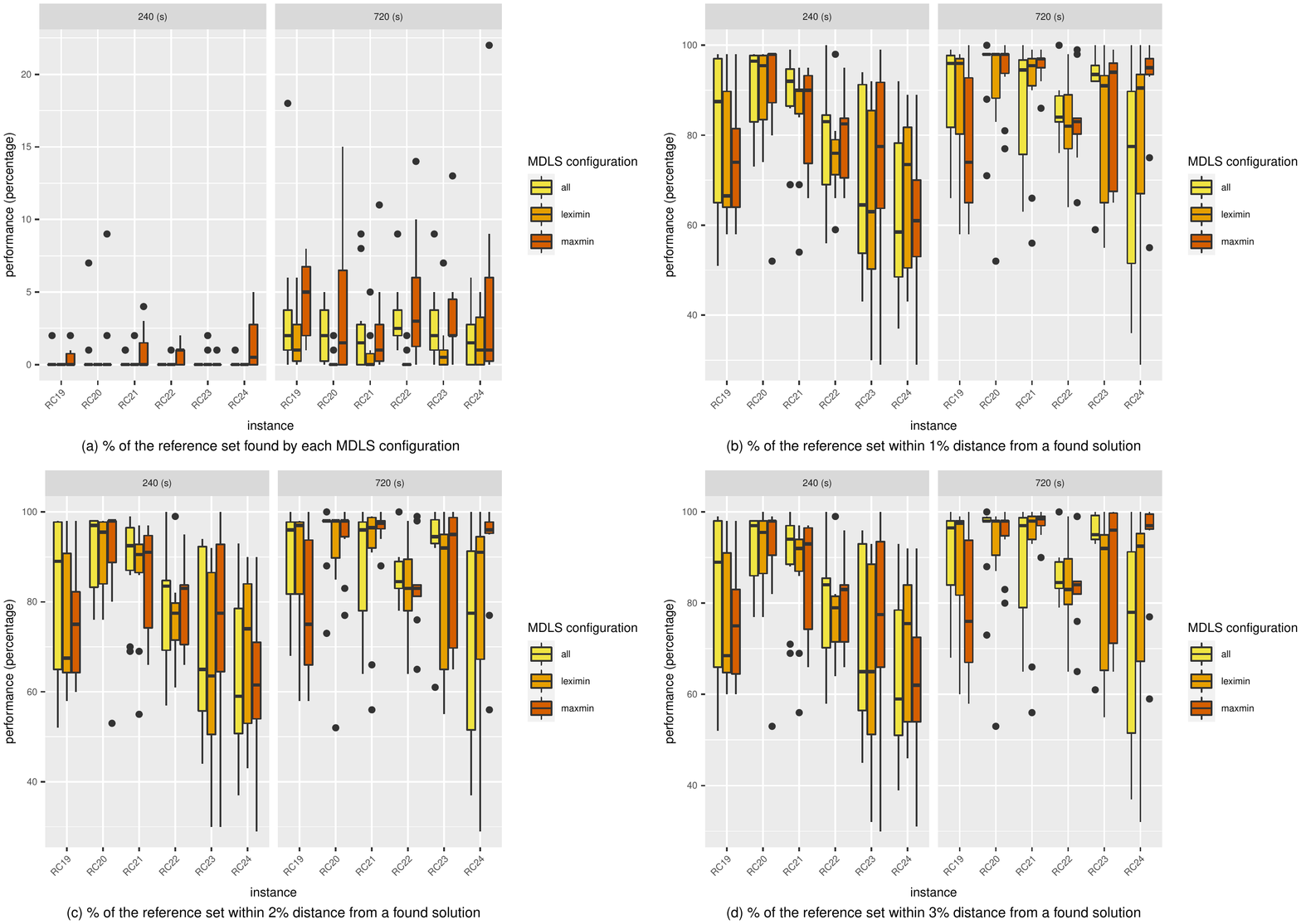}  
  \caption{Proportion of solutions in the reference set found within a 0\%, 1\%, 2\%, 3\% distance (number of sites: 100\_RC, run time: 240 and 720 seconds).}
     \label{RC100}
  \end{figure}
  \end{landscape} 
  
\begin{landscape}

\section{Distribution analysis of maxmimin equivalents}
\label{DistributionAnalysis}  
\vspace{-.5cm}
   \begin{figure}[H]
 \centering
  \includegraphics[width=.80\linewidth]{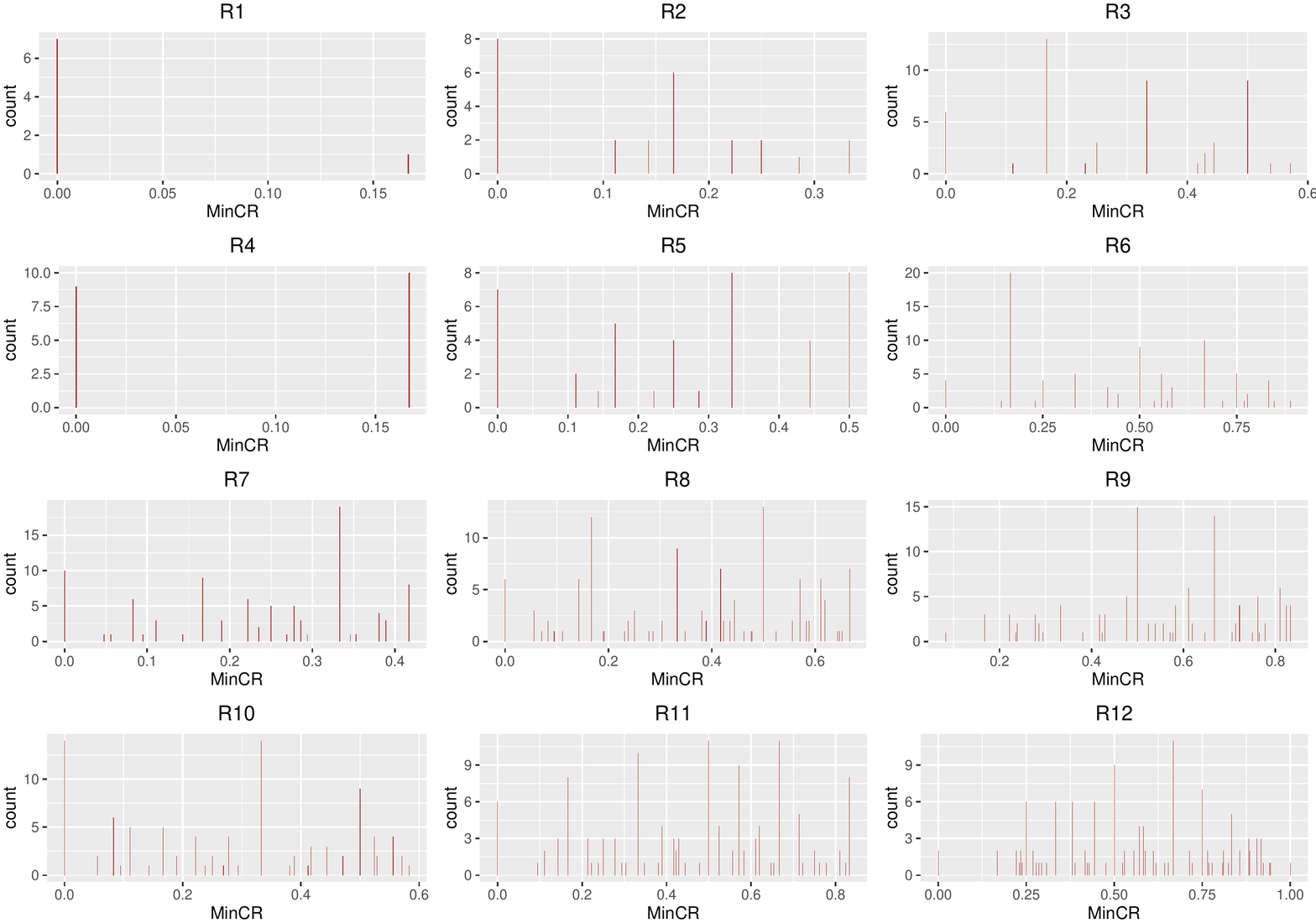}  
  \caption{Distribution analysis for the minimum coverage ratio of non dominated solutions in the leximin- SARP (R1-R12).}
     \label{Rmax-minequalent1}
  \end{figure}
  
   \begin{figure}[H]
 \centering
  \includegraphics[width=0.8\linewidth]{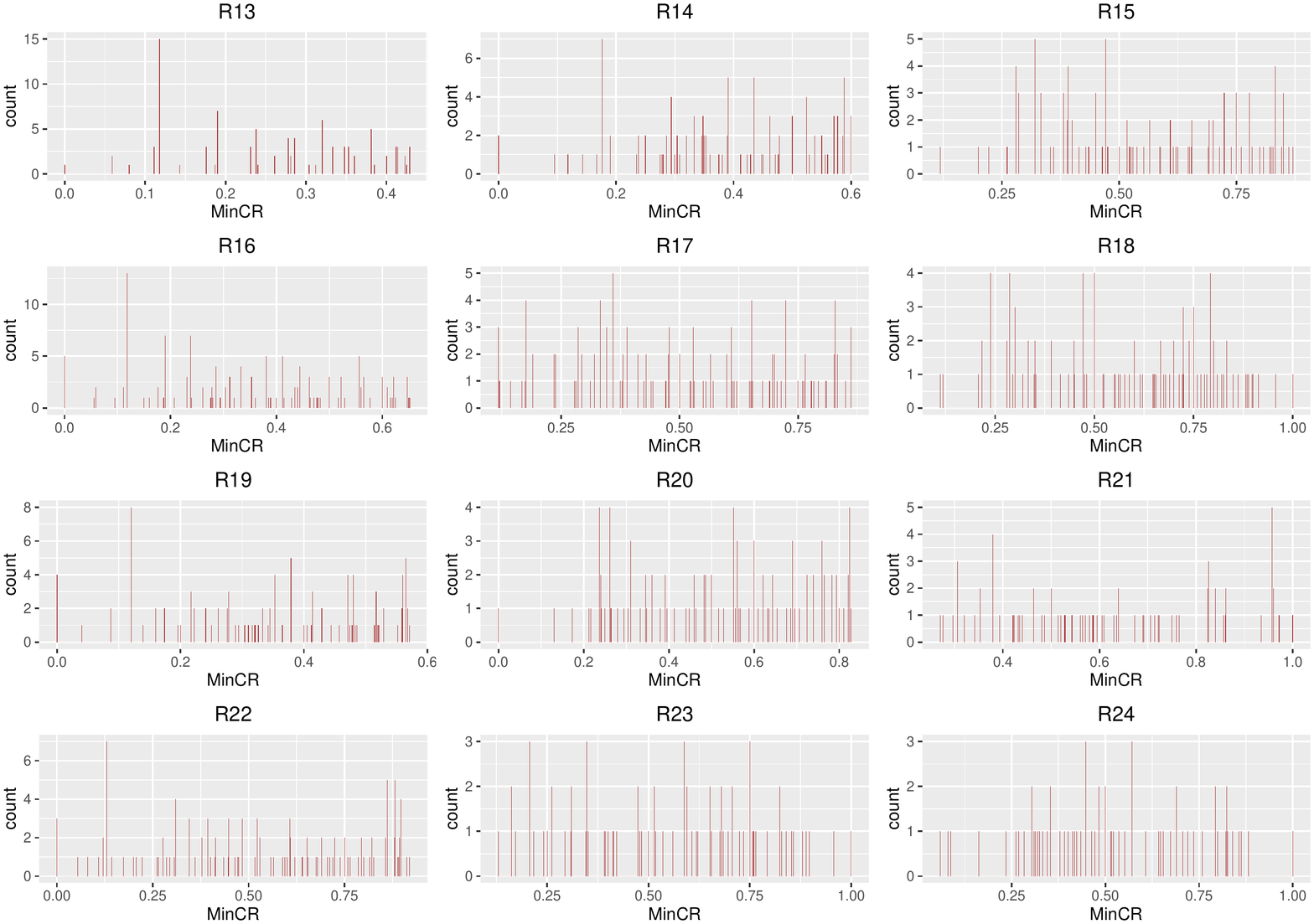}  
  \caption{Distribution analysis for the minimum coverage ratio of non dominated solutions in the leximin- SARP (R13-R24).}
       \label{Rmax-minequalent2}
  \end{figure}
  
     \begin{figure}[H]
 \centering
  \includegraphics[width=0.8\linewidth]{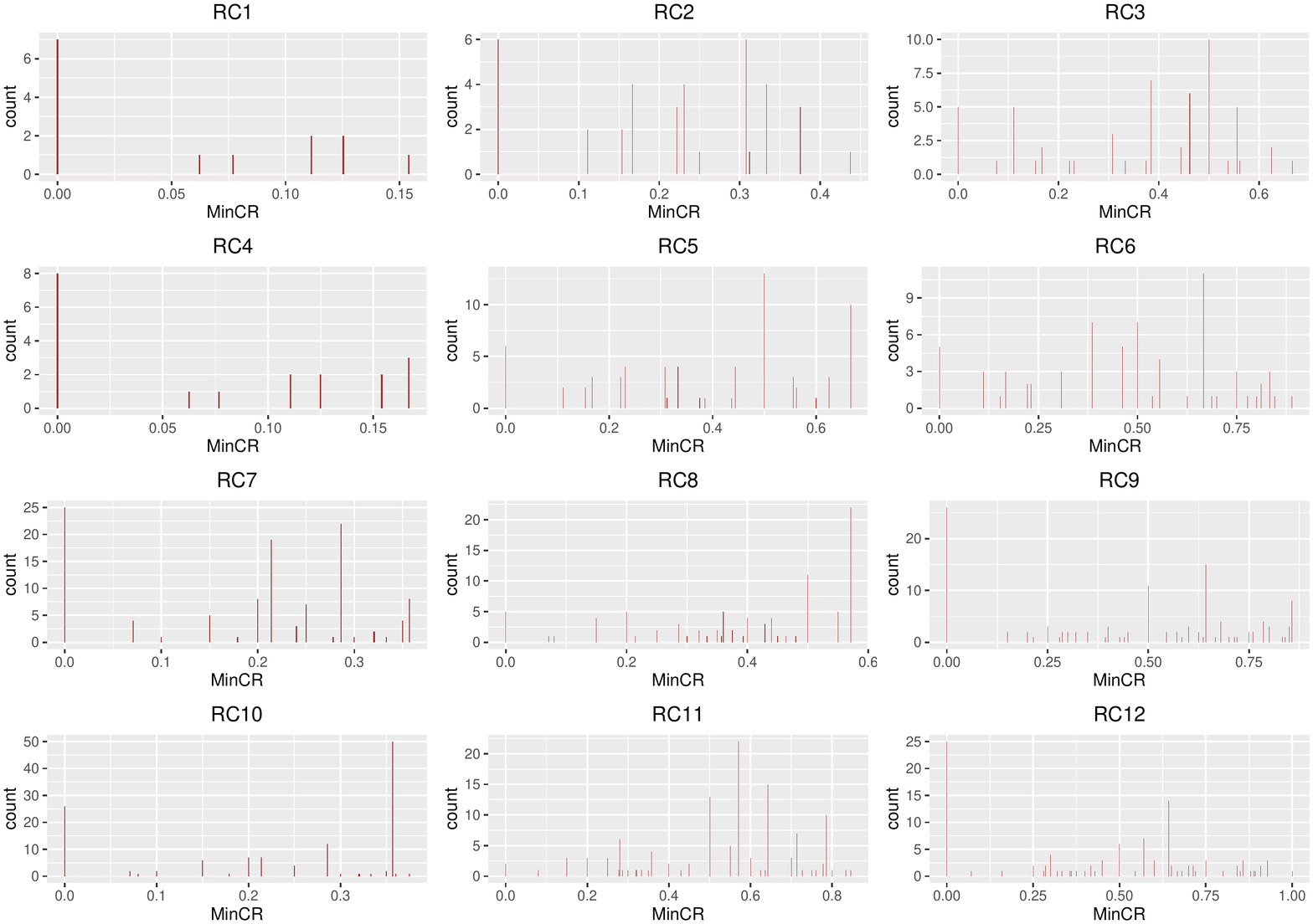}  
  \caption{Distribution analysis for the minimum coverage ratio of non dominated solutions in the leximin- SARP (RC1-RC12).}
       \label{Rmax-minequalent3}
  \end{figure}
     \begin{figure}[H]
 \centering
  \includegraphics[width=0.8\linewidth]{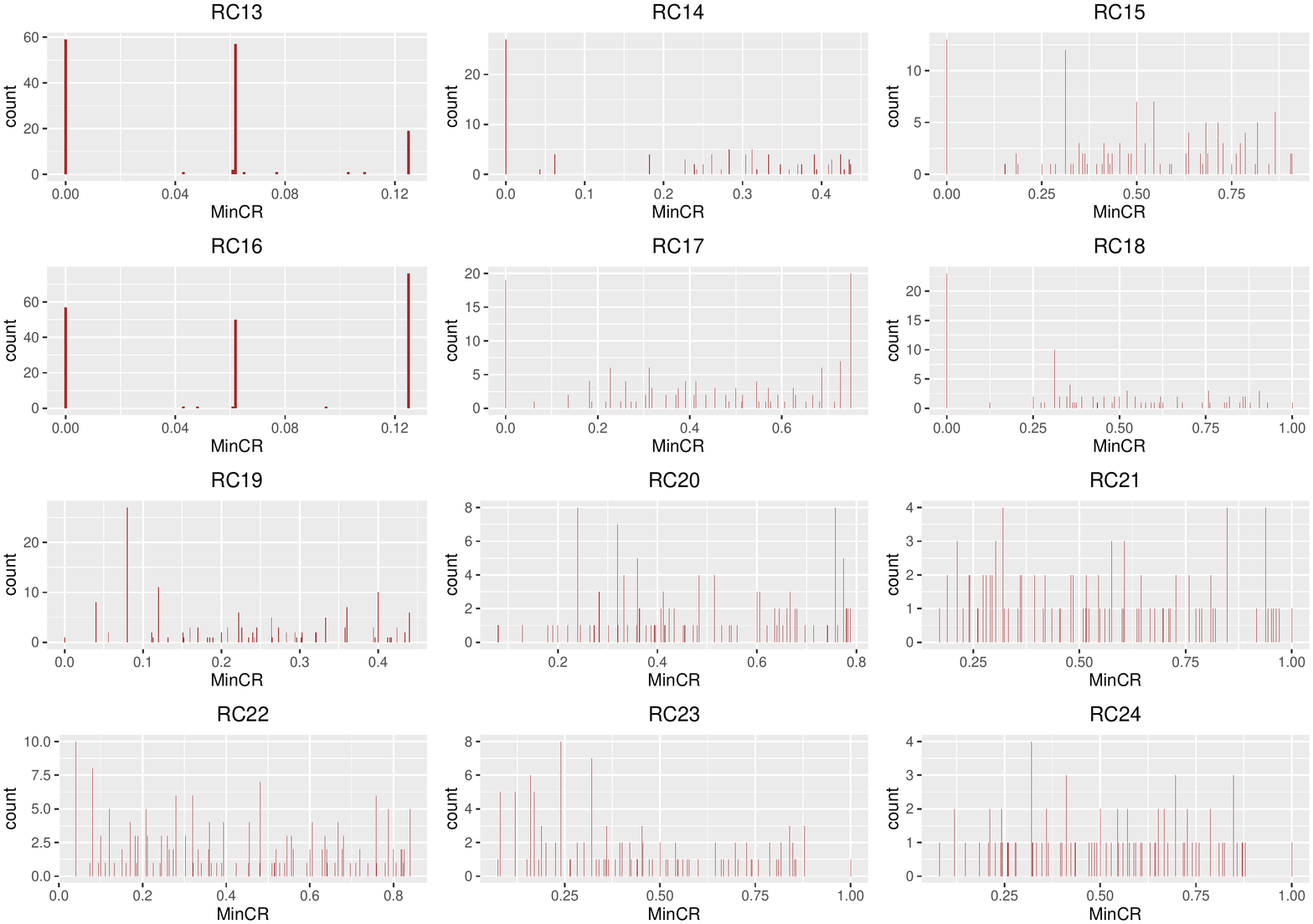}  
  \caption{Distribution analysis for the minimum coverage ratio of non dominated solutions in the leximin- SARP (RC13-RC24).}
       \label{Rmax-minequalent4}
  \end{figure}

\section{Union of Non-dominated solution sets}
\label{Union of Non-dominated}
\vspace{-.5cm}
  \begin{figure}[H]
 \centering
  \includegraphics[width=.8\linewidth]{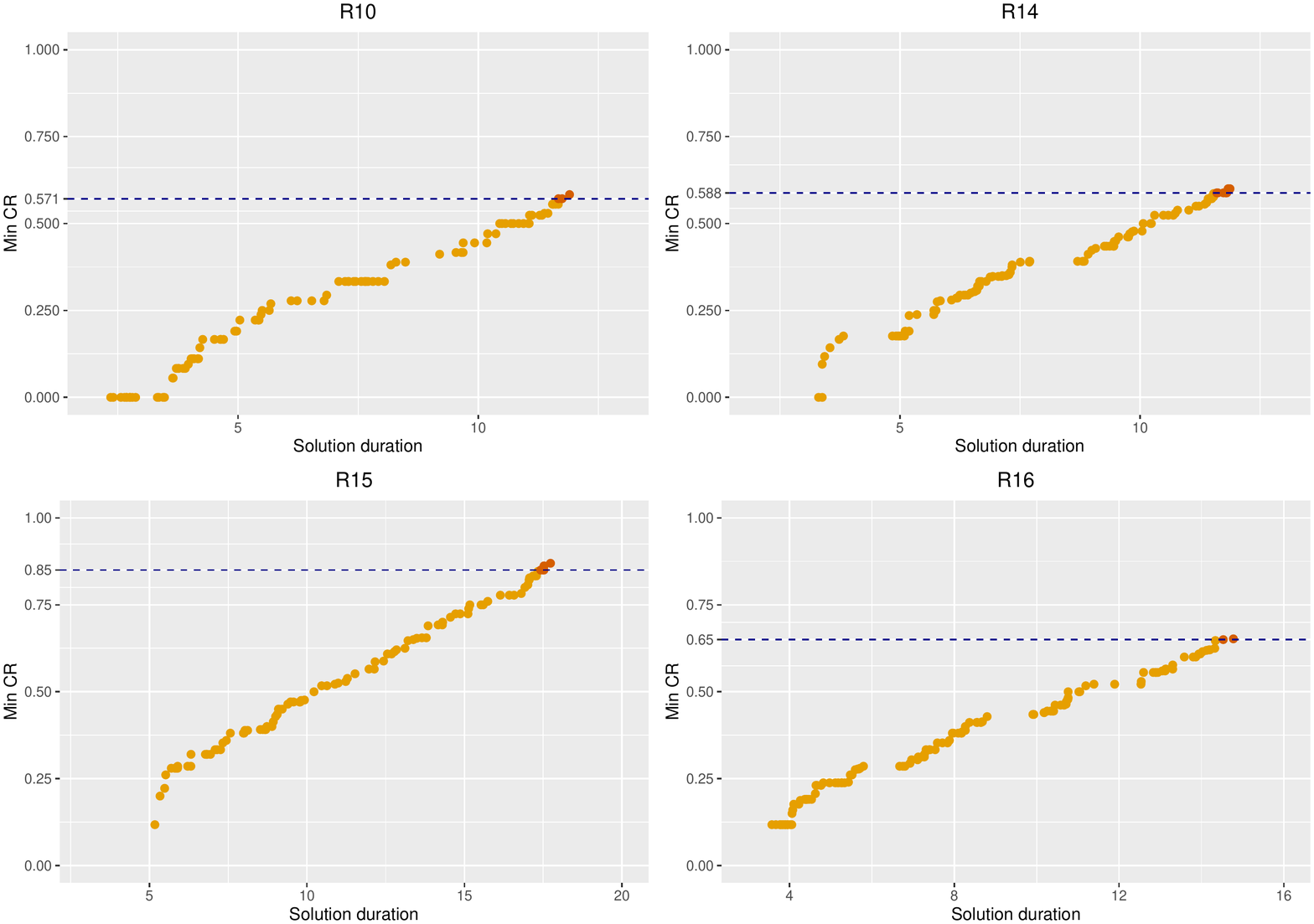}  
  \caption{Union of Non-dominated solution sets of different MDLS variants for the instances R10,R14,R15 and R16}
  \label{RsolutionsDominateBalcik2017}
  \end{figure}
  \begin{figure}[H]
 \centering

  \includegraphics[width=.8\linewidth]{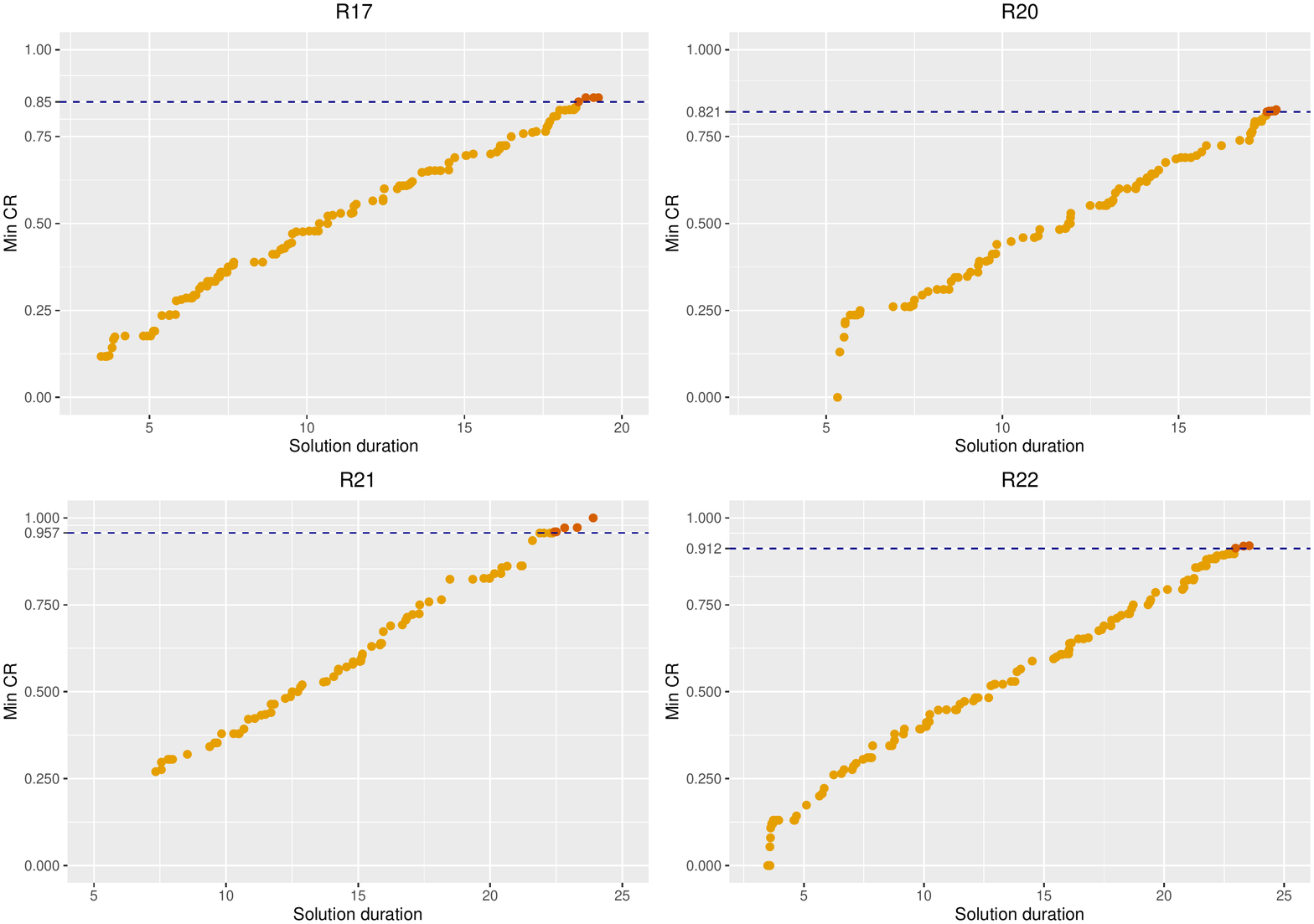}  
  \caption{Union of Non-dominated solution sets of different MDLS variants for the instances R17,R20,R21 and R22}
  \label{RsolutionsDominateBalcik2017-2}
  \end{figure}
  \begin{figure}[H]
 \centering

  \includegraphics[width=.8\linewidth]{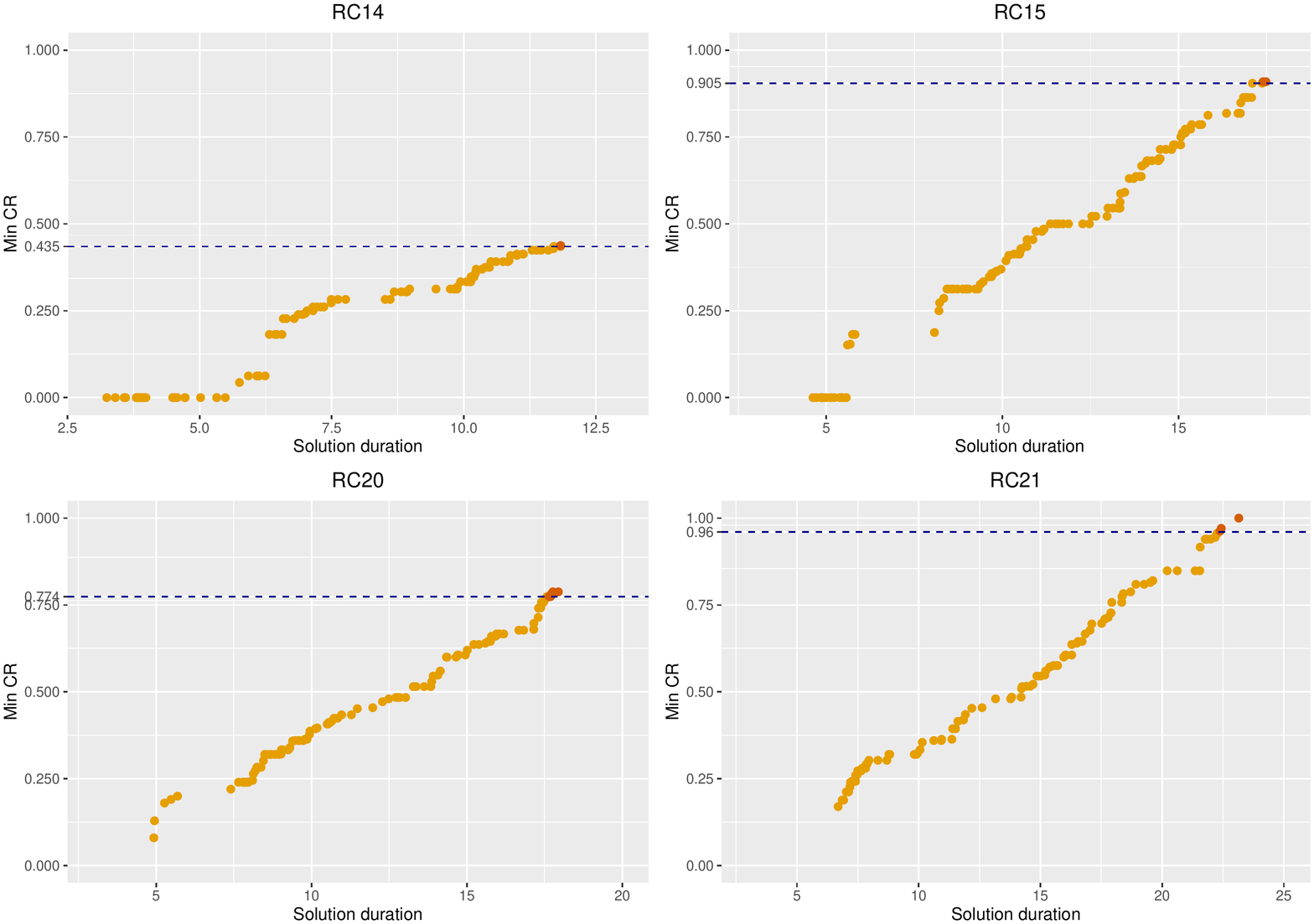}  
  \caption{Union of Non-dominated solution sets of different MDLS variants for the instances RC14,RC15,RC20 and RC21}
  \label{RsolutionsDominateBalcik2017-3}
  \end{figure}

\end{landscape}
\section{Percentage of found best max-min values within individual runs }
  \label{percentageOfBestmax-mins}
\begin{table}[H]
  \centering
  \caption{Percentage of same max-min values reported in \citet{balcik2017site} within 10 runs of each MDLS configuration}
  \scalebox{0.75}{
    \begin{tabular}{|c|cc|cccc|}
    \toprule
    \multirow{2}[4]{*}{Instance} & \multicolumn{1}{p{6.445em}|}{max-min value } & \multicolumn{1}{p{4.165em}|}{CPU time} & \multicolumn{3}{c|}{\% of Runs } & \multicolumn{1}{p{4.165em}|}{CPU time} \\
\cmidrule{4-6}          & \multicolumn{1}{c}{ in Balcik (2017)} & (sec.) & all   & leximin & max-min & (sec.) \\
    \midrule
    R1    & 0.167 & 37    & 100\% & 80\%  & 100\% & \multirow{6}[2]{*}{90} \\
    R2    & 0.333 & 53    & 70\%  & 80\%  & 90\%  &  \\
    R3    & 0.571 & 60    & 20\%  & 0\%   & 20\%  &  \\
    R4    & 0.167 & 49    & 100\% & 100\% & 100\% &  \\
    R5    & 0.5   & 74    & 100\% & 100\% & 100\% &  \\
    R6    & 0.889 & 98    & 100\% & 80\%  & 90\%  &  \\
    \midrule
    R7    & 0.417 & 138   & 90\%  & 100\% & 80\%  & \multirow{6}[2]{*}{180} \\
    R8    & 0.652 & 259   & 90\%  & 40\%  & 90\%  &  \\
    R9    & 0.833 & 466   & 100\% & 60\%  & 100\% &  \\
    R10   & 0.571 & 262   & 60\%  & 20\%  & 80\%  &  \\
    R11   & 0.833 & 451   & 80\%  & 80\%  & 80\%  &  \\
    R12   & 1     & 801   & 100\% & 100\% & 100\% &  \\
    \midrule
    R13   & 0.429 & 317   & 10\%  & 20\%  & 50\%  & \multirow{6}[2]{*}{360} \\
    R14   & 0.588 & 472   & 10\%  & 0\%   & 30\%  &  \\
    R15   & 0.85  & 937   & 80\%  & 50\%  & 90\%  &  \\
    R16   & 0.65  & 554   & 100\% & 30\%  & 100\% &  \\
    R17   & 0.85  & 978   & 70\%  & 60\%  & 70\%  &  \\
    R18   & 1     & 1641  & 100\% & 100\% & 100\% &  \\
    \midrule
    R19   & 0.571 & 1193  & 70\%  & 80\%  & 90\%  & \multirow{6}[2]{*}{720} \\
    R20   & 0.821 & 2154  & 80\%  & 0\%   & 70\%  &  \\
    R21   & 0.957 & 4186  & 100\% & 90\%  & 100\% &  \\
    R22   & 0.912 & 2370  & 80\%  & 0\%   & 70\%  &  \\
    R23   & 1     & 4090  & 100\% & 100\% & 100\% &  \\
    R24   & 1     & 7215  & 100\% & 100\% & 100\% &  \\
    \midrule
    RC1   & 0.154 & 35    & 100\% & 100\% & 100\% & \multirow{6}[2]{*}{90} \\
    RC2   & 0.438 & 51    & 100\% & 100\% & 100\% &  \\
    RC3   & 0.667 & 77    & 100\% & 100\% & 100\% &  \\
    RC4   & 0.167 & 48    & 100\% & 100\% & 100\% &  \\
    RC5   & 0.667 & 101   & 100\% & 100\% & 100\% &  \\
    RC6   & 0.889 & 147   & 100\% & 100\% & 100\% &  \\
    \midrule
    RC7   & 0.357 & 130   & 100\% & 100\% & 100\% & \multirow{6}[2]{*}{180} \\
    RC8   & 0.571 & 285   & 100\% & 100\% & 100\% &  \\
    RC9   & 0.857 & 525   & 70\%  & 40\%  & 70\%  &  \\
    RC10  & 0.357 & 217   & 100\% & 100\% & 100\% &  \\
    RC11  & 0.846 & 403   & 90\%  & 0\%   & 90\%  &  \\
    RC12  & 1     & 794   & 100\% & 100\% & 100\% &  \\
    \midrule
    RC13  & 0.435 & 525   & 100\% & 90\%  & 100\% & \multirow{6}[2]{*}{360} \\
    RC14  & 0.435 & 532   & 100\% & 70\%  & 100\% &  \\
    RC15  & 0.905 & 1103  & 60\%  & 40\%  & 70\%  &  \\
    RC16  & 0.125 & 368   & 100\% & 30\%  & 100\% &  \\
    RC17  & 0.75  & 1036  & 100\% & 100\% & 100\% &  \\
    RC18  & 1     & 1287  & 100\% & 100\% & 100\% &  \\
    \midrule
    RC19  & 0.444 & 648   & 90\%  & 20\%  & 90\%  & \multirow{6}[2]{*}{720} \\
    RC20  & 0.774 & 1593  & 20\%  & 0\%   & 20\%  &  \\
    RC21  & 0.96  & 2634  & 70\%  & 90\%  & 90\%  &  \\
    RC22  & 0.84  & 1857  & 50\%  & 20\%  & 50\%  &  \\
    RC23  & 1     & 2844  & 100\% & 100\% & 100\% &  \\
    RC24  & 1     & 7205  & 100\% & 100\% & 100\% &  \\
    \bottomrule
    \end{tabular}}

\end{table}%

\end{appendices}

\paragraph{Data Availability Statement}The data set that supports the findings of this study is provided by \citet{balcik2017site}.

\pagenumbering{roman}
\small 

\bibliographystyle{apalike}
\bibliography{Ref}


\end{document}